\documentclass{amsart}
\usepackage{amsmath,amssymb,amsthm,bm,bbm}
\usepackage{graphicx,enumerate}
\usepackage{layout}
\usepackage{longtable}
\usepackage[OT2,T1]{fontenc}
\DeclareSymbolFont{cyrletters}{OT2}{wncyr}{m}{n}
\DeclareMathSymbol{\Sha}{\mathalpha}{cyrletters}{"58}
\usepackage[OT2,OT1]{fontenc}
\newcommand\cyr{\renewcommand\rmdefault{wncyr}
\renewcommand\sfdefault{wncyss}
\renewcommand\encodingdefault{OT2}
\normalfont\selectfont}
\DeclareTextFontCommand{\textcyr}{\cyr}

\usepackage{hyperref}

\theoremstyle{plain}
\newtheorem{theorem}{Theorem}[section]
\newtheorem*{theorem-nn}{Theorem}

\newtheorem*{proposition-nn}{Proposition}
\newtheorem{corollary}[theorem]{Corollary}

\newtheorem{conjecture}[theorem]{Conjecture}

\theoremstyle{definition}

\newtheorem{example}[theorem]{Example}

\newtheorem*{acknowledgment}{Acknowledgment}
\theoremstyle{remark}

\newcommand{\bZ}{\mathbbm{Z}}\newcommand{\bQ}{\mathbbm{Q}}
\newcommand{\bG}{\mathbbm{G}}
\newcommand{\bF}{\mathbbm{F}}

\newcommand{\GL}{{\rm GL}}\newcommand{\SL}{{\rm SL}}
\newcommand{\PGL}{{\rm PGL}}\newcommand{\PSL}{{\rm PSL}}
\newcounter{sub}
{\begin{list}{(\arabic{sub})}{\usecounter{sub}%
\setlength{\leftmargin}{2em}}}{\end{list}}

\usepackage{array,booktabs}
\newcommand{\Hind}[2]{\(#1\)\,(#2)}
\newcommand{\Hmult}[3]{\({#1}^{[#3]}\)\,(#2)}
\usepackage{caption}
\title[Rationality problem for norm one tori for small non-Galois extensions]
{Rationality problem for norm one tori for small non-Galois extensions: 
Symmetric, alternating and linear group cases including 
stably rational ${\rm PSL}_2(\mathbbm{F}_{2^d})$ cases with $d=1,2,3$}

\author[A. Hoshi]{Akinari Hoshi}
\address{Department of Mathematics, Niigata University, Niigata 950-2181, Japan}
\email{hoshi@math.sc.niigata-u.ac.jp}

\author[A. Yamasaki]{Aiichi Yamasaki}
\address{Department of Mathematics, Kyoto University, Kyoto 606-8502, Japan}
\email{aiichi.yamasaki@gmail.com}

\thanks{
This work was partially supported by JSPS KAKENHI Grant Numbers 
19K03418, 20H00115, 20K03511, 24K00519, 24K06647.
}


\begin{document}
\begin{abstract}
Let $k$ be a field, $K/k$ be a 
separable field extension of degree $n$ 
and $L/k$ be the Galois closure of $K/k$. 
Let $G={\rm Gal}(L/k)$ and $H={\rm Gal}(L/K)\lneq G$ with $[G:H]=n$. 
Then $H$ is core-free in $G$, 
i.e. $\bigcap_{\sigma\in G} \sigma^{-1}H\sigma=\{1\}$. 
Let $T=R^{(1)}_{K/k}(\bG_m)$ be the norm one torus of $K/k$ 
with ${\rm dim}\, T=n-1$. 
We fix $G$ as an abstract group and let $H$ range 
over all core-free subgroups of $G$ up to conjugacy. 
For each of the following $15$ groups 
$G\simeq S_3$, $S_4$, $S_5$, $S_6$, $A_4$, $A_5$, $A_6$, 
${\rm GL}_2(\mathbbm{F}_q)$ $(q=3,4,5)$, 
${\rm SL}_2(\mathbbm{F}_q)$ $(q=3,5,7)$, 
${\rm PSL}_2(\mathbbm{F}_q)$ $(q=7,8)$, 
by using GAP computations with the aid of PARI/GP, 
we determine, 
for any core-free subgroup $H\leq G$ up to conjugacy, 
whether $T$ is stably $k$-rational, 
not stably but retract $k$-rational, 
or not retract $k$-rational.
The resulting $175$ cases ($15$ are Galois, $160$ are non-Galois) 
split into $7$, $7$, $161$ cases respectively. 
We also give alternative proofs for the known cases 
(i) $H=\{1\}$, i.e. $K/k$ is Galois, and (ii) $[K:k]=[G:H]=n\leq15$ with 
$G\not\simeq 9T27\simeq {\rm PSL}_2(\mathbbm{F}_8)$. 
The $6$ stably $k$-rational cases occur for
$G\simeq S_3\simeq {\rm PSL}_2(\mathbbm{F}_{2})$, 
$A_5\simeq {\rm PSL}_2(\mathbbm{F}_4)$, ${\rm PSL}_2(\mathbbm{F}_8)$. 
Except for the stably $k$-rational case $d=1$ and $H=\{1\}$, 
for $G\simeq {\rm PSL}_2(\mathbbm{F}_{2^d})$ with $d=1,2,3$, 
we prove that $T$ is stably $k$-rational if and only if 
$T$ is retract $k$-rational if and only if 
$(C_2)^d\leq H\leq (C_2)^d\rtimes C_{2^d-1}$, 
and conjecture that this holds for any $d\geq 1$. 
When $k$ is a global field and $T$ is retract $k$-rational, 
e.g. $(G,H)\simeq (\PSL_2(\bF_7),D_4)$, $(\PSL_2(\bF_8),(C_2)^3)$, 
$(\PSL_2(\bF_8),(C_2)^3\rtimes C_7)$, 
the Hasse norm principle holds for $K/k$. 
\end{abstract}

\maketitle

\tableofcontents

\section{Introduction}\label{seInt}

Let $k$ be a field, 
%
$K/k$ be a separable field extension of degree $n$ 
and $L/k$ be the Galois closure of $K/k$. 
Let $G={\rm Gal}(L/k)$ and $H={\rm Gal}(L/K)$ with $[G:H]=n$. 
Then $H$ is {\it core-free in $G$}, 
i.e. $\bigcap_{\sigma\in G} H^\sigma=\{1\}$ 
where $H^\sigma=\sigma^{-1}H\sigma$, that is, 
$H$ contains no normal subgroup of $G$ except for $\{1\}$. 
The Galois group $G$ may be regarded as a transitive subgroup of 
the symmetric group $S_n$ of degree $n$. 
We may assume that 
$H$ is the stabilizer of one of the letters in $G$, 
i.e. $L=k(\theta_1,\ldots,\theta_n)$ and $K=L^H=k(\theta_i)$ 
where $1\leq i\leq n$. 

Let $\overline{k}$ be a fixed separable closure of $k$ and 
$\mathcal{G}={\rm Gal}(\overline{k}/k)$. 
Let $T=R^{(1)}_{K/k}(\bG_m)$ be the norm one torus of $K/k$ 
with ${\rm dim}\, T=n-1$, 
i.e. the kernel of the norm map $R_{K/k}(\bG_m)\rightarrow \bG_m$ where 
$R_{K/k}$ is the Weil restriction 
(see Voskresenskii \cite[page 37, Section 3.12]{Vos98}). 
Let $L/k$ be the Galois closure of $K/k$ with Galois group $G={\rm Gal}(L/k)$. 
Then $T$ splits over $L$, i.e. $T\times_k L\simeq(\bG_{m,L})^{n-1}$. 
Let $\widehat{T}={\rm Hom}(T\times_k \overline{k},\bG_{m,\overline{k}})$ 
be the character lattice of $T$ 
which is a $\mathcal{G}$-lattice, i.e. finitely generated 
$\bZ[\mathcal{G}]$-module which is $\bZ$-free as an abelian group. 
The character lattice $\widehat{T}$ of $T$ 
may be regarded as a $G$-lattice via 
$\widehat{T}\simeq {\rm Hom}(T\times_k L,\bG_{m,L})$ 
with $G\simeq \mathcal{G}/\mathcal{H}$ 
where $\mathcal{H}={\rm Gal}(\overline{k}/L)$ 
and $\widehat{T}\simeq  
J_{G/H}=(I_{G/H})^\circ={\rm Hom}_\bZ(I_{G/H},\bZ)$ 
with $0\to \bZ\to \bZ[G/H]\to J_{G/H}\to 0$ 
which is dual to 
$0\to I_{G/H}\to \bZ[G/H]\xrightarrow{\varepsilon} \bZ\to 0$ 
where $\varepsilon$ is the augmentation map 
(see \cite[Section 4.8]{Vos98}). 
By Endo and Miyata \cite[Theorem 1.6]{EM73} and 
Saltman \cite[Theorem 3.14]{Sal84}, we have
\begin{align*}
T=R^{(1)}_{K/k}(\bG_m)\ {\rm is\ stably\ }k{\text -}{\rm rational}
&\iff [J_{G/H}]^{fl}=0,\\
T=R^{(1)}_{K/k}(\bG_m)\ {\rm is\ retract\ }k{\text -}{\rm rational}
&\iff [J_{G/H}]^{fl}\ {\rm is\ invertible}
\end{align*}
(see Section \ref{sePre} for the flabby class $[J_{G/H}]^{fl}$ of $J_{G/H}$). 
Note that $\bZ[G/H]\simeq \bZ[G/H^\sigma]$ and 
$J_{G/H}\simeq J_{G/H^\sigma}$ for any $\sigma\in G$ 
where $H^\sigma=\sigma^{-1}H\sigma$. 
Hence the stably/retract $k$-rationality problem for $T=R^{(1)}_{K/k}(\bG_m)$ 
is reduced to determining 
the flabby class $[J_{G/H}]^{fl}$ for core-free subgroups $H\leq G$ 
up to conjugacy. 

Let $C_n$ (resp. $D_n$, $A_n$, $S_n$) be the cyclic 
(resp. the dihedral, the alternating, the symmetric) group of 
order $n$ (resp. $2n$, $n!/2$, $n!$).  
Let $V_4\simeq C_2\times C_2$ be the Klein four-group. 
Let $\GL_n(\bF_{p^d})$ (resp. $\SL_n(\bF_{p^d})$, $\PGL_n(\bF_{p^d})$, $\PSL_n(\bF_{p^d})$) be  
the general (resp. special, projective general, projective special) 
linear group of degree $n$ 
over the finite field $\bF_{p^d}$ with $p^d$ elements. 
Let $Z(G)$ be the center of $G$, 
$D(G)$ be the derived (commutator) subgroup of $G$, 
${\rm Syl}_p(G)$ be a $p$-Sylow subgroup of $G$ and 
$N_G(H)$ be the normalizer of $H\leq G$ in $G$. 

The rationality problem for norm one tori is investigated 
by \cite{EM75}, \cite{CTS77}, \cite{Hur84}, \cite{CTS87}, 
\cite{LeB95}, \cite{CK00}, \cite{LL00}, \cite{Flo}, \cite{End11}, 
\cite{HY17}, \cite{HHY20}, \cite{HY21} and \cite{HY24}. 


%
\begin{theorem}[{Endo and Miyata \cite[Theorem 1.5]{EM75}, Saltman \cite[Theorem 3.14]{Sal84}}]\label{th2-2}
Let $L/k$ be a finite Galois field extension and $G={\rm Gal}(L/k)$. 
Then the following conditions are equivalent:\\
{\rm (i)} $R^{(1)}_{L/k}(\bG_m)$ is retract $k$-rational;\\
{\rm (ii)} all the Sylow subgroups of $G$ are cyclic. 
\end{theorem}
\begin{theorem}[{Endo and Miyata \cite[Theorem 2.3]{EM75}, Colliot-Th\'{e}l\`{e}ne and Sansuc \cite[Proposition 3]{CTS77}}]\label{th2-3}
Let $L/k$ be a finite Galois field extension and $G={\rm Gal}(L/k)$. 
Then the following conditions are equivalent:\\
{\rm (i)} $R^{(1)}_{L/k}(\bG_m)$ is stably $k$-rational;\\
{\rm (ii)} $G\simeq C_m$ $(m\geq 1)$ or 
$G\simeq C_m\times \langle x,y\mid x^n=y^{2^d}=1,
yxy^{-1}=x^{-1}\rangle$ $(m\geq 1:odd, n\geq 3:odd, d\geq 1)$ 
with ${\rm gcd}\{m,n\}=1$;\\
{\rm (iii)} $G\simeq \langle s,t\mid s^m=t^{2^d}=1, tst^{-1}=s^r\rangle$ 
$(m\geq 1:odd, d\geq 0)$ with $r^2\equiv 1\pmod{m}$;\\
{\rm (iv)} all the Sylow subgroups of $G$ are cyclic and $H^4(G,\bZ)\simeq \widehat H^0(G,\bZ)$ 
where $\widehat H$ is the Tate cohomology.
\end{theorem}
\begin{theorem}[Endo {\cite[Theorem 2.1]{End11}}]\label{th2-4}
Let $K/k$ be a finite non-Galois separable field extension 
and $L/k$ be the Galois closure of $K/k$. 
Assume that the Galois group of $L/k$ is nilpotent. 
Then the norm one torus $R^{(1)}_{K/k}(\bG_m)$ is not 
retract $k$-rational.
\end{theorem}
\begin{theorem}[Endo {\cite[Theorem 3.1]{End11}}]\label{th2-5}
Let $K/k$ be a finite non-Galois separable field extension 
and $L/k$ be the Galois closure of $K/k$. 
Let $G={\rm Gal}(L/k)$ and $H={\rm Gal}(L/K)\leq G$. 
Assume that all the Sylow subgroups of $G$ are cyclic. 
Then the norm one torus $R^{(1)}_{K/k}(\bG_m)$ is retract $k$-rational, 
and the following conditions are equivalent:\\
{\rm (i)} 
$R^{(1)}_{K/k}(\bG_m)$ is stably $k$-rational;\\
{\rm (ii)} 
$G\simeq C_m\times D_n$ $(m\geq 1:odd,n\geq 3:odd)$ 
with ${\rm gcd}\{m,n\}=1$ and $H\simeq C_2$;\\
{\rm (iii)} 
$H\simeq C_2$ and $G\simeq C_r\rtimes H$ $(r\geq 3:odd)$ 
where $H$ acts non-trivially on $C_r$. 
\end{theorem}

\begin{example}[Corollary of Theorems \ref{th2-2}, \ref{th2-3}, \ref{th2-4}, \ref{th2-5}]\label{ex2-6}
Let $K/k$ be a separable field extension of degree $n\geq 2$ 
and $L/k$ be the Galois closure of $K/k$. 
Then $G={\rm Gal}(L/k)\leq S_n$ transitive, 
$H={\rm Gal}(L/K)\lneq G$ core-free, 
i.e. $\bigcap_{\sigma\in G} H^\sigma=\{1\}$ 
where $H^\sigma=\sigma^{-1}H\sigma$, and $[G:H]=n$. 
Let $T=R^{(1)}_{K/k}(\bG_m)$ be the norm one torus of $K/k$ 
with ${\rm dim}\, T=n-1$.\\
(1) $G\simeq C_n$. 
We have $H=\{1\}$. 
By Theorem \ref{th2-3}, $T$ is stably $k$-rational.\\
(2) $G\simeq D_m$ $(m\geq 3)$. 
We have $H=\{1\}$ or $H\simeq C_2$ with $H\neq Z(D_m)$. 
When $H=\{1\}$, by Theorem \ref{th2-3}, 
$T$ is stably $k$-rational if and only if $m$ is odd. 
When $H\simeq C_2$ with $H\neq Z(D_m)$, by  
Theorem \ref{th2-5} (resp. Theorem \ref{th2-4}),  
if $m$ is odd (resp. $2$-power), then 
$T$ is stably $k$-rational (resp. not retract $k$-rational).\\
(3) $G\simeq F_{pl}=\langle x,y\mid x^p=y^l=1, y^{-1}xy=x^t, 
t=\lambda^{\frac{p-1}{l}}, \langle\lambda\rangle=\bF_p^\times\rangle$ 
where $p$ is a prime and 
$F_{pl}\simeq C_p\rtimes C_l$ $(2<l\mid p-1)$ 
is the Frobenius group of order $pl$. 
We have $H\leq C_l$. 
When $H=\{1\}$ (resp. $H\neq \{1\}$), by 
Theorem \ref{th2-2} and Theorem \ref{th2-3} 
(resp. Theorem \ref{th2-5}),  
$T$ is not stably but retract $k$-rational.\\
(4) $G\simeq Q_{4m}=\langle x,y\mid x^{2m}=y^4=1, x^m=y^2, 
y^{-1}xy=x^{-1}\rangle$ $(m\geq 2)$. 
We have $H=\{1\}$. 
By Theorem \ref{th2-3} (resp. Theorem \ref{th2-2}), 
if $m$ is odd (resp. even), then 
$T$ is stably $k$-rational because 
$G\simeq \langle x^2\rangle\rtimes \langle y\rangle\simeq 
C_m\rtimes C_4$ (resp. not retract $k$-rational 
because $Q_8\leq {\rm Syl}_2(G)$).\\
(5) $G\not\simeq C_m$ and $G$ is nilpotent, 
i.e. the direct product of each $p$-Sylow subgroup of $G$. 
By Theorem \ref{th2-2} (resp. Theorem \ref{th2-4}), 
if $H=\{1\}$ (resp. $H\neq\{1\}$), 
then $T$ is not retract $k$-rational. 
\end{example}

When $G\simeq S_n$ or $A_n$ with $[G:H]=[K:k]=n$, we have:
\begin{theorem}[{Colliot-Th\'{e}l\`{e}ne and Sansuc \cite[Proposition 9.1]{CTS87}, 
\cite[Theorem 3.1]{LeB95}, 
\cite[Proposition 0.2]{CK00}, \cite{LL00}, 
Endo \cite[Theorem 4.1]{End11}, see also 
\cite[Remark 4.2 and Theorem 4.3]{End11}}]\label{thS}
Let $K/k$ be a non-Galois separable field extension 
of degree $n$ and $L/k$ be the Galois closure of $K/k$. 
Assume that ${\rm Gal}(L/k)=S_n$, $n\geq 3$, 
and ${\rm Gal}(L/K)=S_{n-1}$ is the stabilizer of one of the letters 
in $S_n$. 
Then we have:\\
{\rm (i)}\ 
$R^{(1)}_{K/k}(\bG_m)$ is retract $k$-rational 
if and only if $n$ is a prime;\\
{\rm (ii)}\ 
$R^{(1)}_{K/k}(\bG_m)$ is $($stably$)$ $k$-rational 
if and only if $n=3$.
\end{theorem}
\begin{theorem}[Endo {\cite[Theorem 4.4]{End11}, Hoshi and Yamasaki \cite[Corollary 1.11]{HY17}}]\label{thA}
Let $K/k$ be a non-Galois separable field extension 
of degree $n$ and $L/k$ be the Galois closure of $K/k$. 
Assume that ${\rm Gal}(L/k)=A_n$, $n\geq 4$, 
and ${\rm Gal}(L/K)=A_{n-1}$ is the stabilizer of one of the letters 
in $A_n$. 
Then we have:\\
{\rm (i)}\ 
$R^{(1)}_{K/k}(\bG_m)$ is retract $k$-rational 
if and only if $n$ is a prime.\\
{\rm (ii)}\ $R^{(1)}_{K/k}(\bG_m)$ is stably $k$-rational 
if and only if $n=5$.
\end{theorem}

A necessary and sufficient condition for the classification 
of stably/retract $k$-rational norm one tori $T=R^{(1)}_{K/k}(\bG_m)$ 
with $[K:k]=n$ $(n\leq 15$, $n=2^e$ or $n=p$ is a prime$)$ and 
$G={\rm Gal}(L/k)\leq S_n$ 
except for the stable $k$-rationality of $T$ 
with $G\simeq 9T27$ and $G\leq S_p$ for Fermat primes $p\geq 17$ 
was given by Hoshi and Yamasaki \cite{HY21} (for the cases $n=p$ and $n\leq 10$)
and Hasegawa, Hoshi and Yamasaki \cite{HHY20} (for the cases $n=2^e$ and $n=12,14,15$) 
where 
$nTm$ is  the $m$-th transitive subgroup of $S_n$ (see \cite{GAP}). 
%
It seems that the cases where $T=R^{(1)}_{K/k}(\bG_m)$ 
is stably $k$-rational are rare (see Theorem \ref{th2-3}, Theorem \ref{th2-5}, 
Theorem \ref{thS} (ii), Theorem \ref{thA} (ii)). 
We will find two new stably $k$-rational 
$T=R^{(1)}_{K/k}(\bG_m)$ with $G\simeq \PSL_2(\bF_8)$ 
including the remaining case $G\simeq 9T27$ above 
as $(G,H)\simeq (\PSL_2(\bF_8), (C_2)^3\rtimes C_7)$ with $[G:H]=9$ 
in Theorem \ref{mainth}. 


Much of the previous work on the rationality problem for norm one tori 
$T=R^{(1)}_{K/k}(\bG_m)$ 
has been organized by degree: one fixes $n=[K:k]$ and a realization of $G$ 
as a transitive subgroup of $S_n$, and takes
$H={\rm Stab}_G(i)$ to be a stabilizer of one of the letters in $G$. 
In contrast, in the present paper, we fix $G$ as an abstract group and 
let $H$ range over all core-free subgroups of $G$ up to conjugacy. 
In this approach, the extension degree $n=[K:k]=[G:H]$ is not fixed 
in advance but depends on each $H\leq G$.


For each of the $15$ groups listed in Theorem \ref{mainth}, 
we completely determine the stably/retract $k$-rationality 
of the associated norm one torus $T=R^{(1)}_{K/k}(\bG_m)$ 
for every such $H$, 
and it is new except for the known cases (1) $K/k$ is Galois, i.e. $H=1$ 
(see Theorem \ref{th2-2} and Theorem \ref{th2-3}) and 
(2) $n=[K:k]=[G:H]\leq 15$ with 
$G\not\simeq 9T27\simeq {\rm PSL}_2(\mathbbm{F}_8)$ 
as we explained above (although we give an alternative proof of 
these known cases).  
For the definitions of stably/retract $k$-rationalities, 
see Section \ref{sePre}. 
\begin{theorem}\label{mainth}
Let $k$ be a field, $K/k$ be a 
separable field extension of degree $n\geq 2$ 
and $L/k$ be the Galois closure of $K/k$ 
with $G={\rm Gal}(L/k)\leq S_n$ transitive, 
$H={\rm Gal}(L/K)\lneq G$ core-free, 
i.e. $\bigcap_{\sigma\in G} \sigma^{-1}H\sigma=\{1\}$, and $[G:H]=n$. 
Let $T=R^{(1)}_{K/k}(\bG_m)$ be the norm one torus of $K/k$ 
with ${\rm dim}\, T=n-1$. 
Assume that $G$ is isomorphic to one of the following $15$ groups:
\begin{align*}
 S_3,\, S_4,\, S_5,\, S_6,\, 
 A_4,\, A_5,\, A_6,\, 
 {\rm GL}_2(\mathbbm{F}_3),\, 
 {\rm GL}_2(\mathbbm{F}_4),\, 
 {\rm GL}_2(\mathbbm{F}_5),\, 
 {\rm SL}_2(\mathbbm{F}_3),\, 
 {\rm SL}_2(\mathbbm{F}_5),\, 
 {\rm SL}_2(\mathbbm{F}_7),\, 
 {\rm PSL}_2(\mathbbm{F}_7),\, 
 {\rm PSL}_2(\mathbbm{F}_8).
\end{align*}
Then $T$ is 
{\rm (i)} stably $k$-rational, 
{\rm (ii)} not stably but retract $k$-rational 
or 
{\rm (iii)} not retract $k$-rational as specified in Table \ref{Table1}.
\end{theorem}
We note that 
$S_3\simeq {\rm PSL}_2(\mathbbm{F}_2)\simeq {\rm PGL}_2(\mathbbm{F}_2)
\simeq {\rm SL}_2(\mathbbm{F}_2)\simeq {\rm GL}_2(\mathbbm{F}_2)$, 
$A_4\simeq {\rm PSL}_2(\mathbbm{F}_3)$, 
$A_5\simeq {\rm PSL}_2(\mathbbm{F}_5)\simeq {\rm PSL}_2(\mathbbm{F}_4)\simeq 
{\rm PGL}_2(\mathbbm{F}_4) \simeq {\rm SL}_2(\mathbbm{F}_4)$,
$A_6\simeq {\rm PSL}_2(\mathbbm{F}_9)$,
${\rm GL}_2(\mathbbm{F}_4)\simeq A_5\times C_3$, 
${\rm PSL}_2(\mathbbm{F}_7)\simeq {\rm PSL}_3(\mathbbm{F}_2)$, 
${\rm PSL}_2(\mathbbm{F}_8)\simeq {\rm PGL}_2(\mathbbm{F}_8)
\simeq {\rm SL}_2(\mathbbm{F}_8)$.

We also attempted the corresponding computations for 
$S_7$, $A_7$, 
${\rm GL}_2(\mathbbm{F}_7)$, 
${\rm GL}_2(\mathbbm{F}_8)\simeq {\rm PSL}_2(\mathbbm{F}_8)\times C_7$, 
${\rm SL}_2(\mathbbm{F}_{11})$, 
${\rm PSL}_2(\mathbbm{F}_{11})$, 
${\rm PSL}_2(\mathbbm{F}_{16})\simeq {\rm PGL}_2(\mathbbm{F}_{16})
\simeq {\rm SL}_2(\mathbbm{F}_{16})$ 
but were unable to complete them because of limitations in 
memory resources and CPU time.  

In Table \ref{Table1}, we list all core-free subgroups 
$H\lneq G$ up to conjugacy as in Theorem \ref{mainth}
(\# denotes the number of core-free subgroups $H\lneq G$ up to conjugacy), 
$H(n)$ means that $[G:H]=n$, 
$H^{[r]}(n)$ means that there exist $r$ distinct conjugacy classes 
of subgroups isomorphic to $H$ with $[G:H]=n$. 

\begin{longtable}{@{}
>{\raggedright\arraybackslash}p{0.075\textwidth}
>{\centering\arraybackslash}p{0.035\textwidth}
>{\raggedright\arraybackslash}p{0.15\textwidth}
>{\raggedright\arraybackslash}p{0.20\textwidth}
>{\raggedright\arraybackslash}p{0.45\textwidth}@{}}
\caption{Stably/retract $k$-rationality for $T=R^{(1)}_{K/k}(\mathbbm{G}_m)$ with $\widehat{T}\simeq J_{G/H}$ in {Theorem \ref{mainth}}.}
\label{Table1}\\
\toprule
$G$ & \#
& $H$: (i) $T$ is stably $k$-rational
& $H$: (ii) $T$ is not stably but retract $k$-rational
& $H$: (iii) $T$ is not retract $k$-rational\\
\midrule
\endfirsthead

\multicolumn{5}{c}{\tablename\ \thetable\ (continued)}\\
\toprule
$G$ & \#
& $H$: (i) $T$ is stably $k$-rational
& $H$: (ii) $T$ is not stably but retract $k$-rational
& $H$: (iii) $T$ is not retract $k$-rational\\
\midrule
\endhead

\midrule
\endfoot

\bottomrule
\endlastfoot

$S_3$ & 2
 & \Hind{\{1\}}{6}, \Hind{C_2}{3}
 & \multicolumn{1}{c}{--}
 & \multicolumn{1}{c}{--} \\
\addlinespace

$S_4$ & 7
 & \multicolumn{1}{c}{--}
 & \multicolumn{1}{c}{--}
 & \Hind{\{1\}}{24}, \Hmult{C_2}{12}{2},
   \Hind{C_3}{8}, \Hind{V_4}{6}, \Hind{C_4}{6},
   \Hind{S_3}{4}\\
\addlinespace

$S_5$ & 17
 & \multicolumn{1}{c}{--}
 & \Hind{V_4\leq A_5}{30}, \Hind{D_4}{15},
   \Hind{A_4}{10}, \Hind{S_4}{5}
 & \Hind{\{1\}}{120}, \Hmult{C_2}{60}{2},
   \Hind{C_3}{40}, \Hind{V_4\nleq A_5}{30},
   \Hind{C_4}{30}, \Hind{C_5}{24}, \Hind{C_6}{20},
   \Hmult{S_3}{20}{2}, \Hind{D_5}{12},
   \Hind{D_6}{10}, \Hind{C_5\rtimes C_4}{6}\\
\addlinespace

$S_6$ & 54
 & \multicolumn{1}{c}{--}
 & \multicolumn{1}{c}{--}
 & \Hind{\{1\}}{720}, \Hmult{C_2}{360}{3},
   \Hmult{C_3}{240}{2}, \Hmult{V_4}{180}{5},
   \Hmult{C_4}{180}{2}, \Hind{C_5}{144},
   \Hmult{S_3}{120}{4}, \Hmult{C_6}{120}{2},
   \Hmult{(C_2)^3}{90}{2}, \Hind{C_4\times C_2}{90},
   \Hmult{D_4}{90}{4}, \Hind{C_3\times C_3}{80},
   \Hind{D_5}{72}, \Hmult{A_4}{60}{2},
   \Hmult{D_6}{60}{2}, \Hind{C_2\times D_4}{45},
   \Hind{(C_3)^2\rtimes C_2}{40},
   \Hmult{C_3\times S_3}{40}{2},
   \Hind{C_5\rtimes C_4}{36},
   \Hmult{C_2\times A_4}{30}{2},
   \Hmult{S_4}{30}{4},
   \Hmult{(S_3)^2}{20}{2},
   \Hind{(C_3)^2\rtimes C_4}{20},
   \Hmult{C_2\times S_4}{15}{2},
   \Hmult{A_5}{12}{2},
   \Hind{(S_3)^2\rtimes C_2}{10},
   \Hmult{S_5}{6}{2}\\
\addlinespace

$A_4$ & 3
 & \multicolumn{1}{c}{--}
 & \multicolumn{1}{c}{--}
 & \Hind{\{1\}}{12}, \Hind{C_2}{6}, \Hind{C_3}{4}\\
\addlinespace

$A_5$ & 8
 & \Hind{V_4}{15}, \Hind{A_4}{5}
 & \multicolumn{1}{c}{--}
 & \Hind{\{1\}}{60}, \Hind{C_2}{30}, \Hind{C_3}{20},
   \Hind{C_5}{12}, \Hind{S_3}{10}, \Hind{D_5}{6}\\
\addlinespace

$A_6$ & 21
 & \multicolumn{1}{c}{--}
 & \multicolumn{1}{c}{--}
 & \Hind{\{1\}}{360}, \Hind{C_2}{180},
   \Hmult{C_3}{120}{2}, \Hmult{V_4}{90}{2},
   \Hind{C_4}{90}, \Hind{C_5}{72},
   \Hmult{S_3}{60}{2}, \Hind{D_4}{45},
   \Hind{(C_3)^2}{40}, \Hind{D_5}{36},
   \Hmult{A_4}{30}{2},
   \Hind{(C_3)^2\rtimes C_2}{20},
   \Hmult{S_4}{15}{2},
   \Hind{(C_3)^2\rtimes C_4}{10},
   \Hmult{A_5}{6}{2}\\
\addlinespace

$\mathrm{GL}_2(\mathbbm{F}_3)$ & 5
 & \multicolumn{1}{c}{--}
 & \multicolumn{1}{c}{--}
 & \Hind{\{1\}}{48}, \Hind{C_2}{24}, \Hind{C_3}{16},
   \Hmult{S_3}{8}{2}\\
\addlinespace

$\mathrm{GL}_2(\mathbbm{F}_4)$ & 11
 & \Hind{A_4\leq D(G)}{15}
 & \multicolumn{1}{c}{--}
 & \Hind{\{1\}}{180}, \Hind{C_2}{90},
   \Hmult{C_3}{60}{2}, \Hind{V_4}{45}, \Hind{C_5}{36},
   \Hind{S_3}{30}, \Hind{D_5}{18},
   $A_4^{[2]}\nleq D(G)\, (15)$\\
\addlinespace

$\mathrm{GL}_2(\mathbbm{F}_5)$ & 13
 & \multicolumn{1}{c}{--}
 & \multicolumn{1}{c}{--}
 & \Hind{\{1\}}{480}, \Hind{C_2}{240}, \Hind{C_3}{160},
   \Hmult{C_4}{120}{2}, \Hind{C_5}{96}, \Hind{S_3}{80},
   \Hmult{D_5}{48}{2},
   \Hmult{C_5\rtimes C_4}{24}{4}\\
\addlinespace

$\mathrm{SL}_2(\mathbbm{F}_3)$ & 2
 & \multicolumn{1}{c}{--}
 & \multicolumn{1}{c}{--}
 & \Hind{\{1\}}{24}, \Hind{C_3}{8}\\
\addlinespace

$\mathrm{SL}_2(\mathbbm{F}_5)$ & 3
 & \multicolumn{1}{c}{--}
 & \multicolumn{1}{c}{--}
 & \Hind{\{1\}}{120}, \Hind{C_3}{40}, \Hind{C_5}{24}\\
\addlinespace

$\mathrm{SL}_2(\mathbbm{F}_7)$ & 4
 & \multicolumn{1}{c}{--}
 & \multicolumn{1}{c}{--}
 & \Hind{\{1\}}{336}, \Hind{C_3}{112}, \Hind{C_7}{48},
   \Hind{C_7\rtimes C_3}{16}\\
\addlinespace

$\mathrm{PSL}_2(\mathbbm{F}_7)$ & 14
 & \multicolumn{1}{c}{--}
 & \Hind{D_4}{21}, \Hmult{S_4}{7}{2}
 & \Hind{\{1\}}{168}, \Hind{C_2}{84}, \Hind{C_3}{56},
   \Hmult{V_4}{42}{2}, \Hind{C_4}{42}, \Hind{S_3}{28},
   \Hind{C_7}{24}, \Hmult{A_4}{14}{2},
   \Hind{C_7\rtimes C_3}{8}\\
\addlinespace

$\mathrm{PSL}_2(\mathbbm{F}_8)$ & 11
 & \Hind{(C_2)^3}{63},
   \Hind{(C_2)^3\rtimes C_7}{9}
 & \multicolumn{1}{c}{--}
 & \Hind{\{1\}}{504}, \Hind{C_2}{252}, \Hind{C_3}{168},
   \Hind{V_4}{126}, \Hind{S_3}{84}, \Hind{C_7}{72},
   \Hind{C_9}{56}, \Hind{D_7}{36}, \Hind{D_9}{28}\\

\hline
Total
& 175
& \multicolumn{1}{c}{7}
& \multicolumn{1}{c}{7}
& \multicolumn{1}{c}{161}
\end{longtable}

Hence only $14$ of the $175$ cases are retract $k$-rational, 
and exactly $7$ of these are stably $k$-rational.
The $6$ stably $k$-rational cases occur for 
$G\simeq{\rm PSL}_2(\mathbbm{F}_{2^d})$ with $d=1,2,3$ 
because 
$G\simeq S_3\simeq \PSL_2(\bF_2)$, $A_5\simeq \PSL_2(\bF_4)$. 

It is known (due to Ono \cite[page 70]{Ono63}, 
see also Platonov \cite[page 44]{Pla82}, Kunyavskii \cite[Remark 3]{Kun84}, 
Platonov and Rapinchuk \cite[page 307]{PR94}) that 
if $k$ is a global field, i.e. a finite extension of $\bQ$ or $\bF_q(t)$, and 
$T=R^{(1)}_{K/k}(\bG_m)$ is (stably/retract) $k$-rational, 
then the Hasse norm principle holds for $K/k$ 
(see Hoshi, Kanai and Yamasaki \cite{HKY22}, \cite{HKY23}). 
Hence as a corollary of Theorem \ref{mainth} we have: 
\begin{corollary}
If $T=R_{K/k}^{(1)}(\bG_m)$ is retract $k$-rational, 
i.e. the flabby class $[F]=[J_{G/H}]^{fl}$ is invertible 
$($see Theorem \ref{th2-1}$)$, 
as in Theorem \ref{mainth}, e.g. 
\begin{align*}
(G,H)\simeq 
(\PSL_2(\bF_7),D_4),\, 
(\PSL_2(\bF_8),(C_2)^3),\, (\PSL_2(\bF_8),(C_2)^3\rtimes C_7)
\end{align*}
with $[G:H]=21$, $63$, $9$, 
then we get the vanishing 
$H^1(k,{\rm Pic}\, \overline{X})\simeq H^1(G,{\rm Pic}\, X_L)\simeq 
H^1(G,F)\simeq 
\Sha^2_\omega(G,J_{G/H})\simeq {\rm Br}(X)/{\rm Br}(k)\simeq 
{\rm Br}_{\rm nr}(k(X)/k)/{\rm Br}(k)=0$ where 
$X$ is a smooth $k$-compactification of $T$. 
This implies that, 
when $k$ is a global field, 
$A(T)=0$ and $\Sha(T)=0$, i.e. $T$ has the weak approximation property, 
the Hasse principle holds for all torsors $E$ under $T$ and 
the Hasse norm principle holds for $K/k$ 
$($see 
Hoshi, Kanai and Yamasaki \cite[Section 1]{HKY22}, \cite[Section 1]{HKY},  
Hoshi and Yamasaki \cite[Section 4]{HY}$)$. 
\end{corollary}
The groups $H^1(k,{\rm Pic}\, \overline{X})$, $A(T)$ and $\Sha(T)$ 
were investigated by Macedo and Newton \cite{MN22} when $G\simeq A_n$, $S_n$, 
by Hoshi, Kanai and Yamasaki \cite{HKY22}, \cite{HKY23} 
when $[G:H]\leq 15$ and by \cite{HKY} when $G\simeq M_{11}$, $J_1$.\\

More precisely, for the $6$ stably $k$-rational cases with 
$G\simeq S_3\simeq \PSL_2(\bF_2)$, $A_5\simeq \PSL_2(\bF_4)$, $\PSL_2(\bF_8)$ 
as in Theorem \ref{mainth}, 
we prove the following result which implies that 
there exists a rational algebraic $k$-torus $T^\prime=\prod_i R_{k_i/k}(\bG_m)$ (for some $k\subset k_i\subset L$) of dimension $r$ such that 
$T\times T^\prime$ is $k$-rational 
(see Hoshi and Yamasaki \cite[Theorem 1.1]{HY24}, 
see also Section \ref{sePre} for the flabby class $[J_{G/H}]^{fl}$ 
and Theorem \ref{th2-1}): 
\begin{theorem}\label{th1.3}
Let $k$ be a field, $K/k$ be a 
separable field extension of degree $n\geq 2$ 
and $L/k$ be the Galois closure of $K/k$ 
with $G={\rm Gal}(L/k)\leq S_n$ transitive, 
$H={\rm Gal}(L/K)\lneq G$ core-free, 
i.e. $\bigcap_{\sigma\in G} \sigma^{-1}H\sigma=\{1\}$, and $[G:H]=n$. 
Let $T=R^{(1)}_{K/k}(\bG_m)$ be the norm one torus of $K/k$ 
with ${\rm dim}\, T=n-1$ whose function field over $k$ 
is $k(T)\simeq L(J_{G/H})^G$ where 
$J_{G/H}=(I_{G/H})^\circ={\rm Hom}_\bZ(I_{G/H},\bZ)\simeq\widehat{T}\simeq 
{\rm Hom}(T\times_k L,\bG_{m,L})$ 
is the Chevalley module with $0\to \bZ\to \bZ[G/H]\to J_{G/H}\to 0$ 
which is dual to 
$0\to I_{G/H}\to \bZ[G/H]\xrightarrow{\varepsilon} \bZ\to 0$ 
where $\varepsilon$ is the augmentation map.\\
{\rm (1)} When $(G,H)\simeq (S_3, \{1\})\simeq (\PSL_2(\bF_2), \{1\})$ with $[G:H]=6$, 
there exists a flabby $G$-lattice $F$ with 
$[F]=[J_{G/H}]^{fl}$ and ${\rm rank}_\bZ\,F=7$ such that 
an isomorphism of $G$-lattices 
\begin{align*}
\bZ[S_3/C_2]^{\oplus 2}\oplus\bZ[S_3/C_3]\simeq\bZ\oplus F
\end{align*}
holds with rank $r=2\cdot 3+2=1+7=8$.\\
{\rm (2)} When $(G,H)\simeq (S_3, C_2)\simeq (\PSL_2(\bF_2), C_2)$ with $[G:H]=3$, 
there exists a flabby $G$-lattice $F$ with 
$[F]=[J_{G/H}]^{fl}$ and ${\rm rank}_\bZ\,F=4$ such that 
an isomorphism of $G$-lattices 
\begin{align*}
\bZ[S_3/C_2]\oplus\bZ[S_3/C_3]\simeq\bZ\oplus F
\end{align*}
holds with rank $r=3+2=1+4=5$ $($cf. Colliot-Th\'{e}l\`{e}ne and Sansuc \cite[Remarque 1]{CTS77}, 
Hoshi and Yamasaki \cite[Table 8 of Theorem 6.3]{HY17}$)$. 
In particular, we get an isomorphism of $G$-lattices 
\begin{align*}
\bZ[S_3/C_2]^{\oplus 2}\oplus\bZ[S_3/C_3]\simeq\bZ\oplus (\bZ[S_3/C_2]\oplus F)
\end{align*}
with rank $r^\prime=2\cdot 3+2=1+7=8$.\\
{\rm (3)} When $(G,H)\simeq (A_5, V_4)\simeq (\PSL_2(\bF_4),V_4)$ with $[G:H]=15$, 
there exists a flabby $G$-lattice $F$ with 
$[F]=[J_{G/H}]^{fl}$ and ${\rm rank}_\bZ\,F=21$ such that 
an isomorphism of $G$-lattices 
\begin{align*}
\bZ[A_5/C_5]\oplus\bZ[A_5/A_4]^{\oplus 2}\simeq\bZ\oplus F
\end{align*}
holds with rank $r=12+2\cdot 5=1+21=22$.\\ 
{\rm (4)} When $(G,H)\simeq (A_5, A_4)\simeq (\PSL_2(\bF_4),A_4)$ 
with $[G:H]=5$, there exists a flabby $G$-lattice $F$ with 
$[F]=[J_{G/H}]^{fl}$ and ${\rm rank}_\bZ\,F=16$ such that 
an isomorphism of $G$-lattices 
\begin{align*}
\bZ[A_5/C_5]\oplus\bZ[A_5/S_3]\simeq\bZ[A_5/D_5]\oplus F
\end{align*}
holds with rank $r=12+10=6+16=22$.\\
{\rm (5)} When $(G,H)\simeq (\PSL_2(\bF_8), (C_2)^3)$ with $[G:H]=63$, 
there exists a flabby $G$-lattice $F$ with 
$[F]=[J_{G/H}]^{fl}$ and ${\rm rank}_\bZ\,F=73$ such that 
an isomorphism of $G$-lattices 
\begin{align*}
\bZ[G/S_3]\oplus\bZ[G/C_9]\oplus\bZ[G/((C_2)^3\rtimes C_7)]^{\oplus 2}
\simeq\bZ[G/S_3]\oplus\bZ\oplus F
\end{align*}
holds with rank $r=84+56+2\cdot 9=84+1+73=158$.\\
{\rm (6)} When $(G,H)\simeq (\PSL_2(\bF_8), (C_2)^3\rtimes C_7)$ 
with $[G:H]=9$, 
there exists a flabby $G$-lattice $F$ with 
$[F]=[J_{G/H}]^{fl}$ and ${\rm rank}_\bZ\,F=64$ such that 
an isomorphism of $G$-lattices 
\begin{align*}
\bZ[G/C_3]\oplus\bZ[G/C_9]\oplus\bZ[G/D_7]\simeq
\bZ[G/C_3]\oplus\bZ[G/D_9]\oplus F
\end{align*}
holds with rank $r=168+56+36=168+28+64=260$. 

In particular, for the cases $(1)$--$(6)$, $[F]=[J_{G/H}]^{fl}=0$ 
and hence 
$T=R^{(1)}_{K/k}(\bG_m)$ is stably $k$-rational. 
More precisely, there exists a rational algebraic $k$-torus $T^\prime$ 
of dimension $r$ such that $\widehat{T}^\prime$ 
is a permutation $G$-lattice with rank $r$ 
in the left-hand side of the isomorphism, 
i.e. $r=8$, $5$, $22$, $22$, $158$, $260$, 
and $T\times T^\prime$ is $k$-rational. 
\end{theorem}
Motivated by Theorem \ref{mainth} and Theorem \ref{th1.3}, 
we formulate the following conjecture: 
%
\begin{conjecture}\label{con1.4}
Let $k$ be a field, $K/k$ be a 
separable field extension of degree $n\geq 2$ 
and $L/k$ be the Galois closure of $K/k$. 
Then $G={\rm Gal}(L/k)\leq S_n$ transitive, 
$H={\rm Gal}(L/K)\lneq G$ core-free, 
i.e. $\bigcap_{\sigma\in G} H^\sigma=\{1\}$ 
where $H^\sigma=\sigma^{-1}H\sigma$, and $[G:H]=n$. 
Let $T=R^{(1)}_{K/k}(\bG_m)$ be the norm one torus of $K/k$ 
with ${\rm dim}\, T=n-1$. 
When $G\simeq \PSL_2(\bF_{2^d})$ $(d\geq 1)$, 
$T$ is not retract $k$-rational, i.e. $[J_{G/H}]^{fl}$ is not invertible, 
except for the case 
$(G,H)\simeq (S_3,\{1\})$ $(d=1)$ 
and the case 
${\rm Syl}_2(G)\simeq (C_2)^d\leq H\leq 
N_G({\rm Syl}_2(G))\simeq (C_2)^d\rtimes C_{2^d-1}$ 
with $|G|=(2^d+1)2^d(2^d-1)$, 
$[G:H]=(2^d+1)u$, $u\mid 2^d-1$ $(d\geq 1)$. 
For the exceptional cases, $T$ is stably $k$-rational, 
i.e. $[J_{G/H}]^{fl}=0$. 
Moreover, \\ 
{\rm (1)} for $H\simeq  (C_2)^d$ $(d\geq 1)$, 
there exist a flabby $\PSL_2(\bF_{2^d})$-lattice $F$ with 
$[F]=[J_{G/H}]^{fl}$ and ${\rm rank}_\bZ\,F=2^{2d}+2^d+1$ 
and a permutation $\PSL_2(\bF_{2^d})$-lattice $P$ such that 
an isomorphism of $\PSL_2(\bF_{2^d})$-lattices 
\begin{align*}
P\oplus \bZ[G/C_{2^d+1}]\oplus\bZ[G/((C_2)^d\rtimes C_{2^d-1})]^{\oplus 2}\simeq 
P\oplus\bZ\oplus F
\end{align*}
holds with ${\rm rank}_\bZ\,P+2^d(2^d-1)+2\times(2^d+1)={\rm rank}_\bZ\,P+1+(2^{2d}+2^d+1)$;\\ 
{\rm (2)} for $H\simeq (C_2)^d\rtimes C_{2^d-1}$ $(d\geq 1)$, 
there exist a flabby $\PSL_2(\bF_{2^d})$-lattice $F$ with 
$[F]=[J_{G/H}]^{fl}$ and ${\rm rank}_\bZ\, F=2^{2d}$ 
and a permutation $\PSL_2(\bF_{2^d})$-lattice $Q$ such that 
an isomorphism of $\PSL_2(\bF_{2^d})$-lattices 
\begin{align*}
Q\oplus \bZ[G/C_{2^d+1}]\oplus\bZ[G/D_{2^d-1}]\simeq 
Q\oplus\bZ[G/D_{2^d+1}]\oplus F
\end{align*}
holds with ${\rm rank}_\bZ\,Q+2^d(2^d-1)+2^{d-1}(2^d+1)
={\rm rank}_\bZ\,Q+2^{d-1}(2^d-1)+2^{2d}$ where $D_1=C_2$ $(d=1)$.
\end{conjecture}
Note that Theorem \ref{th1.3} proves that 
Conjecture \ref{con1.4} (1) holds for $d=1$, $2$, $3$ with 
${\rm rank}_\bZ\,P=0$, $0$, $84$ and 
Conjecture \ref{con1.4} (2) holds for $d=1$, $2$, $3$ with 
${\rm rank}_\bZ\,Q=0$, $0$, $168$. 

We organize this paper as follows. 
Section \ref{sePre} reviews the necessary background and known results. 
Section \ref{seProof} proves Theorem \ref{mainth} and Theorem \ref{th1.3} and 
describes the new GAP functions used in the proofs. 
Section~\ref{GAPcomp} presents the corresponding GAP computations; 
some of the GAP functions we developed use a union-find algorithm and 
call efficient PARI/GP \cite{PARI2} functions. 
Finally, 
Section \ref{GAPalg} gives the GAP algorithms 
which are also available from \cite{RatProbNorm1Tori} 
(Kyoto University Research Information Repository (KURENAI) 
and the authors' algorithm webpages). 

\begin{acknowledgment}
We thank the referee for a careful reading and valuable suggestions. 
In revising the paper in light of this advice, 
we were led to summarize the main theorem (Theorem \ref{mainth}) 
in Table \ref{Table1}, which has clarified the scope and novelty of 
the theorem and substantially improved the paper. 
\end{acknowledgment}

\section{Preliminaries}\label{sePre}

Let $k$ be a field and $K$ 
be a finitely generated field extension of $k$. 
A field $K$ is called {\it rational over $k$} 
(or {\it $k$-rational} for short) 
if $K$ is purely transcendental over $k$, 
i.e. $K$ is isomorphic to $k(x_1,\ldots,x_n)$, 
the rational function field over $k$ with $n$ variables $x_1,\ldots,x_n$ 
for some integer $n$. 
$K$ is called {\it stably $k$-rational} 
if $K(y_1,\ldots,y_m)$ is $k$-rational for some algebraically 
independent elements $y_1,\ldots,y_m$ over $K$. 
Two fields 
$K$ and $K^\prime$ are called {\it stably $k$-isomorphic} if 
$K(y_1,\ldots,y_m)\simeq K^\prime(z_1,\ldots,z_n)$ over $k$ 
for some algebraically independent elements $y_1,\ldots,y_m$ over $K$ 
and $z_1,\ldots,z_n$ over $K^\prime$. 
When $k$ is an infinite field, 
$K$ is called {\it retract $k$-rational} 
if there is a $k$-algebra $R$ contained in $K$ such that 
(i) $K$ is the quotient field of $R$, and (ii) 
the identity map $1_R : R\rightarrow R$ factors through a localized 
polynomial ring over $k$, i.e. there is an element $f\in k[x_1,\ldots,x_n]$, 
which is the polynomial ring over $k$, and there are $k$-algebra 
homomorphisms $\varphi : R\rightarrow k[x_1,\ldots,x_n][1/f]$ 
and $\psi : k[x_1,\ldots,x_n][1/f]\rightarrow R$ satisfying 
$\psi\circ\varphi=1_R$ (cf. Saltman \cite{Sal84}). 
$K$ is called {\it $k$-unirational} 
if $k\subset K\subset k(x_1,\ldots,x_n)$ for some integer $n$. 
It is not difficult to see that 
``$k$-rational'' $\Rightarrow$ ``stably $k$-rational'' $\Rightarrow$ 
``retract $k$-rational'' $\Rightarrow$ ``$k$-unirational''. 

Let $\overline{k}$ be a fixed separable closure of the base field $k$. 
Let $T$ be an algebraic $k$-torus, 
i.e. a group $k$-scheme with fiber product (base change) 
$T\times_k \overline{k}=
T\times_{{\rm Spec}\, k}\,{\rm Spec}\, \overline{k}
\simeq (\bG_{m,\overline{k}})^n$; 
$k$-form of the split torus $(\bG_m)^n$. 
An algebraic $k$-torus 
$T$ is said to be {\it $k$-rational} (resp. {\it stably $k$-rational}, 
{\it retract $k$-rational}, {\it $k$-unirational}) 
if the function field $k(T)$ of $T$ over $k$ is $k$-rational 
(resp. stably $k$-rational, retract $k$-rational, $k$-unirational). 
Two algebraic $k$-tori $T$ and $T^\prime$ 
are said to be 
{\it birationally $k$-equivalent} 
if their function fields $k(T)$ and $k(T^\prime)$ are 
$k$-isomorphic. 
For an equivalent definition in the language of algebraic geometry, 
see e.g. 
Manin \cite{Man86}, 
Manin and Tsfasman \cite{MT86}, 
Colliot-Th\'{e}l\`{e}ne and Sansuc \cite[Section 1]{CTS07}, 
Beauville \cite{Bea16}, Merkurjev \cite[Section 3]{Mer17}.

Let $L$ be a finite Galois extension of $k$ and $G={\rm Gal}(L/k)$ 
be the Galois group of the extension $L/k$. 
Let $M=\bigoplus_{1\leq i\leq n}\bZ\cdot u_i$ be a $G$-lattice with 
a $\bZ$-basis $\{u_1,\ldots,u_n\}$, 
i.e. finitely generated $\bZ[G]$-module 
which is $\bZ$-free as an abelian group. 
Let $G$ act on the rational function field $L(x_1,\ldots,x_n)$ 
over $L$ with $n$ variables $x_1,\ldots,x_n$ by 
\begin{align}
x_i^{\sigma}=\prod_{j=1}^n x_j^{a_{i,j}},\quad 1\leq i\leq n\label{acts}
\end{align}
for any $\sigma\in G$, when $u_i^{\sigma}=\sum_{j=1}^n a_{i,j} u_j$, 
$a_{i,j}\in\bZ$. 
The field $L(x_1,\ldots,x_n)$ with this action of $G$ will be denoted 
by $L(M)$.
There is the duality between the category of $G$-lattices 
and the category of algebraic $k$-tori which split over $L$ 
(see Ono \cite[Section 1.2]{Ono61}, 
Voskresenskii \cite[page 27, Example 6]{Vos98}, 
Knus, Merkurjev, Rost and Tignol \cite[Proposition 20.17]{KMRT98}). 
In fact, if $T$ is an algebraic $k$-torus which splits over $L$, 
then the character lattice
$\widehat{T}={\rm Hom}(T\times_k \overline{k},\bG_{m,\overline{k}})$ 
of $T$ may be regarded as a $G$-lattice via 
$\widehat{T}\simeq {\rm Hom}(T\times_k L,\bG_{m,L})$ as $G$-lattices. 

Conversely, for a given $G$-lattice $M$, there exists an algebraic 
$k$-torus $T={\rm Spec}(L[M]^G)$ which splits over $L$ 
such that $\widehat{T}\simeq M$ as $G$-lattices. 

The invariant field $L(M)^G$ of $L(M)$ under the action of $G$ 
may be identified with the function field $k(T)$ 
of the algebraic $k$-torus $T$. 
Note that the field $L(M)^G\simeq k(T)$ is always $k$-unirational 
(see \cite[page 40, Example 21]{Vos98}). 
Isomorphism classes of tori of dimension $n$ over $k$ correspond bijectively 
to the elements of the set $H^1(\mathcal{G},\GL_n(\bZ))$ 
where $\mathcal{G}={\rm Gal}(\overline{k}/k)$ since 
${\rm Aut}(\bG_m^n)=\GL_n(\bZ)$. 
The algebraic 
$k$-torus $T$ of dimension $n$ is determined uniquely by the integral 
representation $h : \mathcal{G}\rightarrow \GL_n(\bZ)$ up to conjugacy, 
and the group $h(\mathcal{G})$ is a finite subgroup of $\GL_n(\bZ)$ 
(see \cite[page 57, Section 4.9]{Vos98}). 
The minimal splitting field $L$ of $T$ corresponds to 
$\mathcal{H}$ where 
$\rho: \mathcal{G}\rightarrow G$ and $\mathcal{H}={\rm Ker}(\rho)$ 
with $G\simeq \mathcal{G}/\mathcal{H}$.

A $G$-lattice $M$ is called {\it permutation} $G$-lattice 
if $M$ has a $\bZ$-basis permuted by $G$, 
i.e. $M\simeq \oplus_{1\leq i\leq m}\bZ[G/H_i]$ 
for some subgroups $H_1,\ldots,H_m$ of $G$. 
$M$ is called {\it stably permutation} 
$G$-lattice if $M\oplus P\simeq P^\prime$ 
for some permutation $G$-lattices $P$ and $P^\prime$. 
$M$ is called {\it invertible} 
if it is a direct summand of a permutation $G$-lattice, 
i.e. $P\simeq M\oplus M^\prime$ for some permutation $G$-lattice 
$P$ and a $G$-lattice $M^\prime$. 
We say that two $G$-lattices $M_1$ and $M_2$ are {\it similar} 
if there exist permutation $G$-lattices $P_1$ and $P_2$ such that 
$M_1\oplus P_1\simeq M_2\oplus P_2$. 
The set of similarity classes becomes a commutative monoid 
with respect to the sum $[M_1]+[M_2]:=[M_1\oplus M_2]$ 
and the zero $0=[P]$ where $P$ is a permutation $G$-lattice. 
For a $G$-lattice $M$, there exists a short exact sequence of $G$-lattices 
$0 \rightarrow M \rightarrow P \rightarrow F \rightarrow 0$
where $P$ is permutation and $F$ is flabby which is called a 
{\it flabby resolution} of $M$ 
(see Endo and Miyata \cite[Lemma 1.1]{EM75}, 
Colliot-Th\'el\`ene and Sansuc \cite[Lemma 3]{CTS77}, 
Manin \cite[Appendix, page 286]{Man86}). 
The similarity class $[F]$ of $F$ is determined uniquely 
and is called {\it the flabby class} of $M$. 
We denote the flabby class $[F]$ of $M$ by $[M]^{fl}$. 
We say that $[M]^{fl}$ is invertible if $[M]^{fl}=[E]$ for some 
invertible $G$-lattice $E$. 

For algebraic $k$-tori $T$, we see that 
$[\widehat{T}]^{fl}=[{\rm Pic}\,\overline{X}]$ 
where $X$ is a smooth $k$-compactification of $T$, 
i.e. smooth projective $k$-variety $X$ containing $T$ as a dense 
open subvariety, 
and $\overline{X}=X\times_k\overline{k}$ 
(see Voskresenskii \cite[Section 4, page 1213]{Vos69}, \cite[Section 3, page 7]{Vos70}, 
\cite{Vos74}, \cite[Section 4.6]{Vos98}, Kunyavskii \cite[Theorem 1.9]{Kun07} and 
Colliot-Th\'el\`ene \cite[Theorem 5.1, page 19]{CT07} for any field $k$).

The flabby class $[M]^{fl}=[\widehat{T}]^{fl}$ 
plays a crucial role in the rationality problem 
for $L(M)^G\simeq k(T)$ 
as follows (see also 
Colliot-Th\'el\`ene and Sansuc \cite[Section 2]{CTS77}, 
\cite[Proposition 7.4]{CTS87}, 
Voskresenskii \cite[Section 4.6]{Vos98}, 
Kunyavskii \cite[Theorem 1.7]{Kun07}, 
Colliot-Th\'el\`ene \cite[Theorem 5.4]{CT07}, 
Hoshi and Yamasaki \cite[Section 1]{HY17}): 
%
\begin{theorem}
\label{th2-1}
Let $L/k$ be a finite Galois extension with Galois group $G={\rm Gal}(L/k)$ 
and $M$ and $M^\prime$ be $G$-lattices. 
Let $T$ and $T^\prime$ be algebraic $k$-tori with $\widehat{T}\simeq M$ 
and $\widehat{T}^\prime\simeq M^\prime$, 
i.e. $L(M)^G\simeq k(T)$ and $L(M^\prime)^G\simeq k(T^\prime)$.\\
{\rm (i)} $(${\rm Endo and Miyata} \cite[Theorem 1.6]{EM73}$)$ 
$[M]^{fl}=0$ if and only if $L(M)^G$ is stably $k$-rational.\\
{\rm (ii)} $(${\rm Voskresenskii} \cite[Theorem 2]{Vos74}$)$ 
$[M]^{fl}=[M^\prime]^{fl}$ if and only if $L(M)^G$ and $L(M^\prime)^G$ 
are stably $k$-isomorphic. 
{\rm (iii)} $(${\rm Saltman} \cite[Theorem 3.14]{Sal84}$)$ 
$[M]^{fl}$ is invertible if and only if $L(M)^G$ is 
retract $k$-rational.
\end{theorem}

Let $K/k$ be a separable field extension of degree $n$ 
and $L/k$ be the Galois closure of $K/k$. 
Let $G={\rm Gal}(L/k)$ and $H={\rm Gal}(L/K)$ with $[G:H]=n$. 
The norm one torus $R^{(1)}_{K/k}(\bG_m)$ splits over $L$, 
and has the Chevalley module $J_{G/H}$ as its character lattice 
and the field $L(J_{G/H})^G$ as its function field 
where $J_{G/H}=(I_{G/H})^\circ={\rm Hom}_\bZ(I_{G/H},\bZ)$ 
is the dual lattice of $I_{G/H}={\rm Ker}\ \varepsilon$ and 
$\varepsilon : \bZ[G/H]\rightarrow \bZ$ is the augmentation map 
(see \cite[Section 4.8]{Vos98}). 
We have the exact sequence $0\rightarrow \bZ\rightarrow \bZ[G/H]
\rightarrow J_{G/H}\rightarrow 0$ and ${\rm rank}_\bZ\,J_{G/H}=n-1$. 
Write $J_{G/H}=\oplus_{1\leq i\leq n-1}\bZ x_i$. 
Then the action of $G$ on $L(J_{G/H})=L(x_1,\ldots,x_{n-1})$ is 
of the form (\ref{acts}). 

\section{Proofs of Theorem {\ref{mainth}} and Theorem {\ref{th1.3}}}\label{seProof}

Let $K/k$ be a 
separable field extension of degree $n$ 
and $L/k$ be the Galois closure of $K/k$. 
Let $T=R^{(1)}_{K/k}(\bG_m)$ be the norm one torus of $K/k$. 
Assume that $G={\rm Gal}(L/k)$ and $H={\rm Gal}(L/K)\lneq G$. 
Then $H$ is core-free in $G$, i.e. 
$\bigcap_{\sigma\in G} H^\sigma=\{1\}$ 
where $H^\sigma=\sigma^{-1}H\sigma$, that is,  
$H$ contains no normal subgroup of $G$ except for $\{1\}$. 

We may assume that $H$ is the stabilizer of one of the letters in $G$, 
i.e. $L=k(\theta_1,\ldots,\theta_n)$ and $K=L^H=k(\theta_i)$ 
where $1\leq i\leq n$. 

We made the following new GAP \cite{GAP} functions 
(see Section \ref{GAPalg} for details of the GAP codes) 
in order to prove Theorem \ref{mainth} and Theorem \ref{th1.3} 
by investigating the Chevalley module $J_{G/H}$ 
and by applying the function 
{\tt IsInvertibleF} (see Hoshi and Yamasaki 
\cite[Section 5.2, Algorithm F2]{HY17}) to determine whether the flabby class 
$[F]=[J_{G/H}]^{fl}$ is invertible or not, i.e. 
$T=R_{K/k}^{(1)}(\bG_m)$ is retract $k$-rational or not 
(see Theorem \ref{th2-1}). 
Moreover, by a method given in \cite[Chapter 5 and Chapter 8]{HY17} 
with the aid of the new function {\tt StablyPermutationFCheckPFromBasePari} 
which uses a union-find algorithm and requires PARI/GP \cite{PARI2},  
we can confirm that 
$[F]=0$, i.e. $T=R_{K/k}^{(1)}(\bG_m)$ is stably $k$-rational 
(see the proof of Theorem \ref{mainth} below).\\ 


{\tt Hcandidates(}$G${\tt )} returns subgroups $H$ of $G$ which satisfy 
$\bigcap_{\sigma\in G} H^\sigma=\{1\}$ where $H^\sigma=\sigma^{-1}H\sigma$ 
(hence $H$ contains no normal subgroup of $G$ except for $\{1\}$).\\

{\tt Norm1TorusJTransitiveGroup(}$d,n${\tt )} returns the Chevalley module $J_{G/H}$ for the $n$-th transitive 
subgroup $G = {}_dT_n \leq S_d$ of degree $d$ where $H$ is the stabilizer of one of the 
letters in $G$. (The input and output of this function is the same as the function {\tt Norm1TorusJ(}$d,n${\tt )}  
which is given in Hoshi and Yamasaki \cite[Chapter 8]{HY17} but this function is more efficient.)\\ 

{\tt Norm1TorusJCoset(}$G,H${\tt )} returns the Chevalley module $J_{G/H}$ 
for a group $G$ and a subgroup $H\leq G$.\\ 

{\tt StablyPermutationFCheckPPari(}$G,L_1,L_2${\tt )} returns the same as {\tt StablyPermutationFCheckP(}$G,L_1,L_2${\tt )} 
but using efficient PARI/GP functions (e.g. matker, matsnf) \cite{PARI2}. 
(See \cite[Section 5.7, Method III]{HY17}. 
This function applies a union-find algorithm and it also requires PARI/GP \cite{PARI2}.)\\

{\tt StablyPermutationFCheckPFromBasePari(}$G,m_i,L_1,L_2${\tt )} 
returns the same as\\
{\tt StablyPermutationFCheckPPari(}$G,L_1,L_2${\tt )} 
but with respect to $m_i=\mathcal{P}^\circ$ instead of the original $\mathcal{P}^\circ$ 
as in Hoshi and Yamasaki \cite[Equation (4) in Section 5.1]{HY17}. 
(See \cite[Section 5.7, Method III]{HY17}. 
This function applies a union-find algorithm and it also requires PARI/GP \cite{PARI2}.)\\

Some related GAP algorithms are also available 
from \cite{RatProbNorm1Tori} (Kyoto University Research Information Repository (KURENAI) and the authors' algorithm webpages).\\ 

{\it Proof of Theorem \ref{mainth}.}\\

For $G\simeq S_3\simeq\PSL_2(\bF_2)\simeq \PGL_2(\bF_2)\simeq 
\SL_2(\bF_2)\simeq \GL_2(\bF_2)$ with $|G|=6$, 
there exist $4$ subgroups 
$H_1=\{1\}, H_2\simeq C_2, H_3\simeq C_3, H_4\simeq S_3\leq G$ 
up to conjugacy. 
By the condition $\bigcap_{\sigma\in G} H_i^\sigma=\{1\}$, 
we should take core-free $H_i$ $(i=1,2)$ (i.e. $H_i\not\simeq C_3, S_3)$. 
By applying {\tt IsInvertibleF}, we find that 
$[J_{S_3/H_i}]^{fl}$ $(i=1,2)$ is invertible as an $S_3$-lattice. 
Moreover, it follows from Theorem \ref{th2-3} and Theorem \ref{th2-5} that 
$T$ is stably $k$-rational for $i=1,2$. 
Alternatively, 
by a method given in \cite[Chapter 5 and Chapter 8]{HY17}, 
we get an isomorphism of $S_3$-lattices for $i=1$ $(H_1=\{1\})$: 
\begin{align*}
\bZ[S_3/C_2]^{\oplus 2}\oplus\bZ[S_3/C_3]\simeq\bZ\oplus F
\end{align*}
with rank $2\cdot 3+2=1+7=8$ 
where $[F]=[J_{S_3/\{1\}}]^{fl}$ with ${\rm rank}_\bZ\, F=7$. 
Similarly, for $i=2$ $(H_2\simeq C_2)$, we have: 
\begin{align*}
\bZ[S_3/C_2]\oplus\bZ[S_3/C_3]\simeq\bZ\oplus F
\end{align*}
with rank $3+2=1+4=5$ where 
$[F]=[J_{S_3/C_2}]^{fl}$ with ${\rm rank}_\bZ\, F=4$. 
This implies that $[F]=[J_{S_3/H_i}]^{fl}=0$ and hence 
$T$ is stably $k$-rational for $i=1,2$ 
(see Section \ref{GAPcomp} for GAP computations).\\

For $G\simeq S_4\simeq\PGL_2(\bF_3)$ with $|G|=24$ 
(resp. $G\simeq A_4\simeq \PSL_2(\bF_3)$ with $|G|=12$), 
there exist $11$ (resp. $5$) 
subgroups $H\leq G$ up to conjugacy. 
Out of $11$ (resp. $5$), we get 
$s=7$ core-free subgroups $H_1=\{1\},\ldots,H_{7}\simeq S_3\leq G$
(resp. $s=3$ core-free subgroups $H_1=\{1\}$, $H_2\simeq C_2$, $H_3\simeq C_3\leq G$)
with the condition $\bigcap_{\sigma\in G} H_i^\sigma=\{1\}$ $(1\leq i\leq s)$. 
By applying the function {\tt Norm1TorusJCoset(}$G$,$H_i${\tt )}, 
we obtain $J_{G/H_i}$ in GAP. 
Then by the function {\tt IsInvertibleF} (see Hoshi and Yamasaki 
\cite[Section 5.2, Algorithm F2]{HY17}) 
we find that $[J_{G/H_i}]^{fl}$ is not invertible for each $H_i\lneq G$ $(1\leq i\leq s)$. 
It follows from Theorem \ref{th2-1} that $T$ is not retract $k$-rational 
(see Section \ref{GAPcomp} for GAP computations).\\

For $G\simeq S_5\simeq\PGL_2(\bF_5)$ with $|G|=120$, 
there exist $19$ subgroups 
$H_1=\{1\},\ldots,H_{18}\simeq A_5,H_{19}=G\leq G$ up to conjugacy. 
By the condition $\bigcap_{\sigma\in G} H_i^\sigma=\{1\}$, 
we should take core-free $H_i$ $(1\leq i\leq 17)$ (i.e. $H_i\not\simeq A_5, S_5)$. 
By applying {\tt IsInvertibleF}, we find that $[J_{S_5/H_i}]^{fl}$ $(1\leq i\leq 17)$ 
is not invertible except for $i=5,12,14,17$. 
This implies that $T$ is not retract $k$-rational 
except for $i=5,12,14,17$. 
Note that $H_5\simeq H_6\simeq V_4$ and 
$H_5\leq D(G)\simeq A_5$ although $H_6\cap D(G)\simeq C_2$ 
where $D(G)$ is the derived (commutator) subgroup of $G$. 
For the cases $i=5,12,14,17$, we see that $[J_{S_5/H_i}]^{fl}$ is invertible
and hence $T$ is retract $k$-rational. 
We find that $G=30T25$ and $H_5\simeq V_4$ with $[G:H_5]=30$, 
$G=15T10$ and $H_{12}\simeq D_4$ with $[G:H_{12}]=15$, 
$G=10T12$ and $H_{14}\simeq A_4$ with $[G:H_{14}]=10$,  
$G=5T5$ and  $H_{17}\simeq S_4$ with $[G:H_{17}]=5$. 
Hence it follows from \cite[Theorem 1.2 (1), (4)]{HHY20} 
and \cite[Theorem 4.1]{End11} 
(see also \cite[Theorem 1.10]{HY17}) that 
$T$ is not stably but retract $k$-rational for $i=12,14,17$. 
Here we give a proof of not only the remaining case $i=5$ but also $i=12,14,17$ together. 
By \cite[Algorithm 4.1 (2)]{HHY20}, we obtain that 
$[F_i]=[J_{S_5/H_i}]^{fl}\neq 0$  
with ${\rm rank}_\bZ\, F_i=151, 76, 41, 16$ for $i=5,12,14,17$ respectively. 
This implies that $T$ is not stably but retract $k$-rational for $i=5,12,14,17$ 
(see Section \ref{GAPcomp} for GAP computations).\\

For $G\simeq S_6$ with $|G|=720$ 
(resp. $G\simeq A_6\simeq\PSL_2(\bF_9)$ with $|G|=360$), 
there exist $56$ (resp. $22$) subgroups 
$H_1=\{1\},\ldots,H_{55}\simeq A_6,H_{56}=G\leq G$
(resp. $H_1=\{1\},\ldots,H_{22}=G\leq G$) up to conjugacy. 
By the condition $\bigcap_{\sigma\in G} H_i^\sigma=\{1\}$, 
we should take core-free $H_i$ $(1\leq i\leq s)$ with $s=54$ (resp. $s=21$). 
By {\tt IsInvertibleF}, we find that 
$[J_{G/H_i}|_{{\rm Syl_2}(G)}]^{fl}$ $(1\leq i\leq s)$ is not invertible as a ${\rm Syl}_2(G)$-lattice 
except for $i=27, 34, 43, 44, 48, 49$ (resp. $i=11, 17, 18$). 
For the cases $i=27, 34, 43, 44, 48, 49$ (resp. $i=11, 17, 18$), 
we obtain that 
$[J_{G/H_i}|_{{\rm Syl_3}(G)}]^{fl}$ is not invertible as a ${\rm Syl}_3(G)$-lattice. 
This implies that $[J_{G/H_i}]^{fl}$ is not invertible as a $G$-lattice for each $1\leq i\leq s$ 
(see Hoshi and Yamasaki \cite[Lemma 2.17]{HY17}, 
see also Colliot-Th\'{e}l\`{e}ne and Sansuc \cite[Remarque R2, page 180]{CTS77}). 
It follows from Theorem \ref{th2-1} that $T$ is not retract $k$-rational 
(see Section \ref{GAPcomp} for GAP computations).\\

For $G\simeq A_5\simeq\PSL_2(\bF_5)\simeq \PSL_2(\bF_4)\simeq 
\PGL_2(\bF_4)\simeq \SL_2(\bF_4)$ with $|G|=60$, there exist $9$ subgroups 
$H_1=\{1\},\ldots,H_{9}=G\leq G$ up to conjugacy. 
As in the case of $G\simeq A_4$, 
by applying {\tt IsInvertibleF}, we find that 
$[J_{A_5/H_i}]^{fl}$ $(1\leq i\leq 8)$ is not invertible as an $A_5$-lattice 
except for $i=4, 8$. 
For $i=4, 8$, we have $G=15T5$ and $H_4\simeq V_4$ with $[G:H_4]=15$, 
and $G=5T4$ and $H_8\simeq A_4$ with $[G:H_8]=5$ 
where 
$nTm$ is  the $m$-th transitive subgroup of $S_n$ (see \cite{GAP}). 
Hence it follows from \cite[Theorem 1.2 (4)]{HHY20} and \cite[Theorem 1.10]{HY17} that 
$T$ is stably $k$-rational for $i=4,8$. 
Alternatively, 
by a method given in \cite[Chapter 5 and Chapter 8]{HY17}, 
we get an isomorphism of $A_5$-lattices for $i=4$: 
\begin{align*}
\bZ[A_5/H_5]\oplus\bZ[A_5/H_8]^{\oplus 2}\simeq\bZ\oplus F
\end{align*}
with rank $12+2\cdot 5=1+21=22$ 
where $H_5\simeq C_5$, $H_8\simeq A_4$, 
$[F]=[J_{A_5/H_4}]^{fl}$ with ${\rm rank}_\bZ\, F=21$. 
Similarly, for $i=8$, we have: 
\begin{align*}
\bZ[A_5/H_5]\oplus\bZ[A_5/H_6]\simeq\bZ[A_5/H_7]\oplus F
\end{align*}
with rank $12+10=6+16=22$ 
where $H_5\simeq C_5$, $H_6\simeq S_3$, $H_7\simeq D_5$, 
$[F]=[J_{A_5/H_8}]^{fl}$ with ${\rm rank}_\bZ\, F=16$. 
This implies that $[F]=[J_{A_5/H_i}]^{fl}=0$ and hence 
$T$ is stably $k$-rational for $i=4,8$ 
(see Section \ref{GAPcomp} for GAP computations).\\

For $G\simeq \GL_2(\bF_3)$ with $|G|=48$ 
(resp. $G\simeq \SL_2(\bF_3)$ with $|G|=24$, 
$G\simeq \SL_2(\bF_5)$ with $|G|=120$), 
there exist $16$ (resp. $7$, $12$) subgroups $H\leq G$ up to conjugacy. 
Out of $16$ (resp. $7$, $12$), we get $s=5$ (resp. $2$, $3$) core-free  
subgroups 
$H_1=\{1\},\ldots,H_{5}\simeq S_3\leq G$ 
(resp. $H_1=\{1\}, H_{2}\simeq C_3\leq G$, 
$H_1=\{1\},\ldots,H_{3}\simeq C_5\leq G$)  
with the condition $\bigcap_{\sigma\in G} H_i^\sigma=\{1\}$ $(1\leq i\leq s)$. 
Applying {\tt IsInvertibleF}, we find that $[J_{G/H_i}]^{fl}$ $(1\leq i\leq s)$ is not invertible. 
It follows from Theorem \ref{th2-1} that $T$ is not retract $k$-rational 
(see Section \ref{GAPcomp} for GAP computations).\\

For $G\simeq \GL_2(\bF_4)\simeq A_5\times C_3$ with $|G|=180$, 
there exist $21$ subgroups $H\leq G$ up to conjugacy. 
Out of $21$, we get $11$ core-free subgroups 
$H_1=\{1\},\ldots,H_{11}\simeq A_4\leq G$ 
with the condition $\bigcap_{\sigma\in G} H_i^\sigma=\{1\}$ $(1\leq i\leq 11)$. 
Applying {\tt IsInvertibleF}, we find that $[J_{G/H_i}|_{{\rm Syl_2}(G)}]^{fl}$ $(1\leq i\leq 11)$ is not invertible 
except for $i=5,11$. 
For $i=5$, we have $H_5\simeq V_4$ and $[J_{G/H_5}|_{{\rm Syl_3}(G)}]^{fl}$ 
is not invertible. 
This implies that $T$ is not retract $k$-rational except for the case $i=11$. 
For the case $i=11$, we have $H_{11}\simeq A_4$. 
Note that $H_9\simeq H_{10}\simeq H_{11}\simeq A_4$ and 
$D(G)\cap H_9\simeq D(G)\cap H_{10}\simeq V_4$ although $H_{11}\leq D(G)\simeq A_5$
where $D(G)$ is the derived (commutator) subgroup of $G$. 
We also see $G=15T16$ and $[G:H_{11}]=15$. 
Hence it follows from \cite[Theorem 1.2 (4)]{HHY20} that $T$ is stably $k$-rational. 
Alternatively, by \cite[Algorithm 4.1 (3)]{HHY20}, we get 
$[F]=[J_{G/H_{11}}]^{fl}$ with ${\rm rank}_\bZ\,F=36$, 
$[F^\prime]=[F]^{fl}$ with ${\rm rank}_\bZ\,F^\prime=14$ 
and 
$[F^{\prime\prime}]=[F^\prime]^{fl}$ with $F^{\prime\prime}\simeq \bZ$
and hence $[F]=[F^\prime]=[F^{\prime\prime}]=0$ 
(see Section \ref{GAPcomp} for GAP computations).\\

For $G\simeq \SL_2(\bF_7)$ with $|G|=336$ 
(resp. $G\simeq \GL_2(\bF_5)$ with $|G|=480$),  
there exist $19$ (resp. $48$) subgroups $H\leq G$ up to conjugacy. 
Out of $19$ (resp. $48$), we get $s=4$ (resp. $13$) core-free subgroups 
$H_1=\{1\},\ldots,H_4\simeq C_7\rtimes C_3\leq G$ 
(resp. $H_1=\{1\},\ldots,H_{13}\simeq C_5\rtimes C_4\leq G$) 
with the condition $\bigcap_{\sigma\in G} H_i^\sigma=\{1\}$ $(1\leq i\leq s)$. 
Applying {\tt IsInvertibleF}, 
we find that 
$[J_{G/H_i}|_{{\rm Syl_2}(G)}]^{fl}$ $(1\leq i\leq s)$ is not invertible as a ${\rm Syl}_2(G)$-lattice. 
This implies that $[J_{G/H_i}]^{fl}$ is not invertible as a $G$-lattice for each $1\leq i\leq s$ 
(see Hoshi and Yamasaki \cite[Lemma 2.17]{HY17}, 
see also Colliot-Th\'{e}l\`{e}ne and Sansuc \cite[Remarque R2, page 180]{CTS77}). 
It follows from Theorem \ref{th2-1} that $T$ is not retract $k$-rational 
(see Section \ref{GAPcomp} for GAP computations).\\ 


For $G\simeq \PSL_2(\bF_7)\simeq \PSL_3(\bF_2)$ with $|G|=168$, 
there exist $15$ subgroups 
$H_1=\{1\},\ldots,H_{15}=G\leq G$ up to conjugacy. 
By the condition $\bigcap_{\sigma\in G} H_i^\sigma=\{1\}$, 
we should take core-free $H_i$ $(1\leq i\leq 14)$. 
By {\tt IsInvertibleF}, we find that $[J_{G/H_i}]^{fl}$ $(1\leq i\leq 14)$ is not invertible 
except for $i=9, 13, 14$.
For the cases $i=9,13,14$, we see that $[J_{G/H_i}]^{fl}$ is invertible. 
We also find that $G=21T14$ and $H_9\simeq D_4$ with $[G:H_9]=21$, 
$G=7T5$ and $H_{13}\simeq H_{14}\simeq S_4$ with $[G:H_{13}]=[G:H_{14}]=7$. 
Hence it follows from \cite[Theorem 8.5 (ii)]{HY17} that 
$T$ is not stably but retract $k$-rational for $i=13,14$. 
Here we give a proof of not only the remaining case $i=9$ but also $i=13,14$ together. 
By \cite[Algorithm 4.1 (2)]{HHY20}, we obtain that $[F_i]=[J_{G/H_i}]^{fl}\neq 0$  
with ${\rm rank}_\bZ\, F_i=148, 36, 36$ for $i=9,13,14$ respectively. 
This implies that $T$ is not stably but retract $k$-rational for $i=9,13,14$ 
(see Section \ref{GAPcomp} for GAP computations).\\ 

For $G\simeq \PSL_2(\bF_8)$ with $|G|=504$, 
there exist $12$ subgroups 
$H_1=\{1\},\ldots,H_{12}=G\leq G$ up to conjugacy. 
By the condition $\bigcap_{\sigma\in G} H_i^\sigma=\{1\}$, 
we should take core-free $H_i$ $(1\leq i\leq 11)$. 
By {\tt IsInvertibleF}, we find that 
$[J_{G/H_i}|_{{\rm Syl_2}(G)}]^{fl}$ $(1\leq i\leq 11)$ is not invertible as a ${\rm Syl}_2(G)$-lattice 
except for $i=7, 11$. 
For the cases $i=7$, 
$G\leq S_{63}$ is transitive and $H_7\simeq (C_2)^3$ with $[G:H_7]=63$. 
For $i=11$, we have 
$G=9T27$ and $H_{11}\simeq (C_2)^3\rtimes C_7$ with $[G:H_{11}]=9$ 
where 
$nTm$ is  the $m$-th transitive subgroup of $S_n$ (see \cite{GAP}). 
By a method given in \cite[Chapter 5 and Chapter 8]{HY17} 
with the aid of the function {\tt StablyPermutationFCheckPFromBasePari} 
which uses a union-find algorithm and requires PARI/GP \cite{PARI2} 
(see the third paragraph of Section \ref{seProof}), 
we get an isomorphism of $\PSL_2(\bF_8)$-lattices for $i=7$: 
\begin{align*}
\bZ[G/H_5]\oplus\bZ[G/H_8]\oplus\bZ[G/H_{11}]^{\oplus 2}
\simeq\bZ[G/H_5]\oplus\bZ\oplus F
\end{align*}
with rank $84+56+2\cdot 9=84+1+73=158$
where $H_5\simeq S_3$, $H_8\simeq C_9$, 
$H_{11}\simeq (C_2)^3\rtimes C_7$, 
$[F]=[J_{G/H_7}]^{fl}$ with ${\rm rank}_\bZ\,F=73$. 
Similarly, for $i=11$, we have: 
\begin{align*}
\bZ[G/H_3]\oplus\bZ[G/H_8]\oplus\bZ[G/H_9]\simeq
\bZ[G/H_3]\oplus\bZ[G/H_{10}]\oplus F
\end{align*}
with rank $168+56+36=168+28+64=260$ 
where $H_3\simeq C_3$, $H_8\simeq C_9$, 
$H_9\simeq D_7$, $H_{10}\simeq D_9$, 
$[F]=[J_{G/H_{11}}]^{fl}$ with ${\rm rank}_\bZ\,F=64$. 
This implies that $[F]=[J_{G/H_i}]^{fl}=0$ and hence 
$T$ is stably $k$-rational for $i=7,11$ 
(see Section \ref{GAPcomp} for GAP computations).
\qed\\

{\it Proof of Theorem \ref{th1.3}.}\\ 

The first half of the theorem follows from the proof of Theorem \ref{mainth}. 
We will prove the last statement. 
For a flabby resolution 
\begin{align*}
0\to J_{G/H}\to P \to F\to 0
\end{align*}
of $J_{G/H}$ where $P$ is a permutation $G$-lattice and 
the flabby class $[F]=[J_{G/H}]^{fl}$ of $J_{G/H}$, 
we get a modified flabby resolution of $J_{G/H}$: 
\begin{align*}
0\to J_{G/H}\to Q\oplus P\to Q\oplus F\to 0
\end{align*}
where $Q\oplus P$ and $Q\oplus F$ are permutation $G$-lattices 
given in Theorem \ref{th1.3} with ${\rm rank}_\bZ\,(Q\oplus F)=r$. 
In particular, $Q\oplus F$ is an invertible $G$-lattice. 
Hence it follows from 
Lenstra \cite[Proposition 1.5]{Len74} 
(see also Ono \cite[Proposition 1.2.2]{Ono63}, 
Endo and Miyata \cite[Proposition 1.10]{EM73}, 
Swan \cite[Lemma 3.1]{Swa10}) that 
$L(J_{G/H}\oplus Q\oplus F)^G\simeq L(Q\oplus P)^G$. 
This implies that 
$T\times T^\prime$ and $T^{\prime\prime}$ 
are birationally $k$-equivalent and hence $T\times T^\prime$ is $k$-rational 
where $T^\prime$ and $T^{\prime\prime}$ are $k$-rational algebraic $k$-tori 
whose function fields over $k$ are 
$k(T^\prime)\simeq L(Q\oplus F)^G$ and 
$k(T^{\prime\prime})\simeq L(Q\oplus P)^G$ respectively 
(see Hoshi, Kang and Kitayama \cite[Proof of Theorem 6.5]{HKK14}). 
\qed 
\section{GAP computations}\label{GAPcomp}

The following computations were done using GAP, Version 4.9.3 \cite{GAP}. 
The GAP algorithms and related ones are available 
as {\tt FlabbyResolutionFromBase.gap} in \cite{RatProbNorm1Tori}. 

\begin{verbatim}
gap> Read("FlabbyResolutionFromBase.gap");

gap> S3:=SymmetricGroup(3); # G=S3=PSL(2,2)=PGL(2,2)=SL(2,2)=GL(2,2) with |G|=6
Sym( [ 1 .. 3 ] )
gap> S3cs:=ConjugacyClassesSubgroups2(S3); # subgroups H of G up to conjugacy 
[ Group( () )^G, Group( [ (2,3) ] )^G, Group( [ (1,2,3) ] )^G, 
  Group( [ (1,2,3), (2,3) ] )^G ]
gap> S3H:=Hcandidates(S3); # exclude not core-free H=C3, S3
[ Group(()), Group([ (2,3) ]) ]
gap> S3J:=List(S3H,x->Norm1TorusJCoset(S3,x));;
gap> List(S3J,IsInvertibleF); # T is retract k-rational
[ true, true ]

gap> mi:=SearchCoflabbyResolutionBaseLowRank(TransposedMatrixGroup(S3J[1]),0);; # H=1
gap> F:=FlabbyResolutionFromBase(S3J[1],mi).actionF;;
gap> Rank(F.1); # [F]=[J_{G/H1}]^{fl} with rank F=7
7
gap> ll:=PossibilityOfStablyPermutationFFromBase(S3J[1],mi);
[ [ 0, 2, 1, -1, -1 ] ]
gap> l:=ll[1]; # possibility for 2Z[G/C2]+Z[G/C3]=Z+F with rank 2x3+2=1+7=8
[ 0, 2, 1, -1, -1 ]
gap> bp:=StablyPermutationFCheckPFromBase(S3J[1],mi,Nlist(l),Plist(l));;
gap> Length(bp);
14
gap> Length(bp[1]);
8
gap> SearchPRowBlocks(bp);
rec( bpBlocks := [ [ 1, 2, 3 ], [ 4, 5, 6, 7, 8, 9, 10, 11, 12, 13, 14 ] ], 
  rowBlocks := [ [ 1 ], [ 2, 3, 4, 5, 6, 7, 8 ] ] )
gap> r1:=SearchPFilterRowBlocks(bp,[1..3],[1],2);;
gap> Length(r1);
6
gap> r2:=SearchPFilterRowBlocks(bp,[4..14],[2..8],3);;
gap> Length(r2);
12
gap> P:=SearchPMergeRowBlock(r1,r2);
[ [ 1, 1, 1, 0, 0, 0, 1, 1 ], 
  [ 1, 0, 0, 0, 1, 0, 0, 0 ], 
  [ 0, 1, 0, 0, -1, 0, 1, 0 ], 
  [ 1, 0, 0, 0, 0, 1, 0, 0 ], 
  [ 0, 0, 1, 0, 0, -1, 0, 1 ], 
  [ 0, 1, 0, 0, -1, 0, 0, 1 ], 
  [ 0, 0, 1, 0, 1, 0, 0, 0 ], 
  [ 0, 0, 1, 1, 0, 0, 0, 0 ] ]
gap> Determinant(P);
1
gap> StablyPermutationFCheckMatFromBase(S3J[1],mi,Nlist(l),Plist(l),P); 
# T is stably k-rational
true

gap> mi:=SearchCoflabbyResolutionBaseLowRank(TransposedMatrixGroup(S3J[2]),0);; # H=C2
gap> F:=FlabbyResolutionFromBase(S3J[2],mi).actionF;;
gap> Rank(F.1); # [F]=[J_{G/H2}]^{fl} with rank F=4
4
gap> ll:=PossibilityOfStablyPermutationFFromBase(S3J[2],mi);
[ [ 0, 1, 1, -1, -1 ] ]
gap> l:=ll[1]; # possibility for Z[G/C2]+Z[G/C3]=Z+F with rank 3+2=1+4=5 
[ 0, 1, 1, -1, -1 ]
gap> bp:=StablyPermutationFCheckPFromBase(S3J[2],mi,Nlist(l),Plist(l));;
gap> Length(bp);
6
gap> Length(bp[1]);
5
gap> SearchPRowBlocks(bp);
rec( bpBlocks := [ [ 1, 2 ], [ 3, 4, 5, 6 ] ], 
  rowBlocks := [ [ 1 ], [ 2, 3, 4, 5 ] ] )
gap> r1:=SearchPFilterRowBlocks(bp,[1..2],[1],2);;
gap> Length(r1);
3
gap> r2:=SearchPFilterRowBlocks(bp,[3..6],[2..5],2);;
gap> Length(r2);
4
gap> P:=SearchPMergeRowBlock(r1,r2);
[ [ 1, 1, 1, 1, 1 ], 
  [ 1, 1, 0, 1, 0 ], 
  [ 1, 0, 1, 0, 1 ], 
  [ 0, 1, 1, 0, 1 ], 
  [ 0, 1, 1, 1, 0 ] ]
gap> Determinant(P);
-1
gap> StablyPermutationFCheckMatFromBase(S3J[2],mi,Nlist(l),Plist(l),P);
# T is stably k-rational
true

gap> List(S3H,StructureDescription); # for checking core-free H
[ "1", "C2" ]

gap> S4:=SymmetricGroup(4); # G=S4=PGL(2,3) with |G|=24
Sym( [ 1 .. 4 ] )
gap> S4cs:=ConjugacyClassesSubgroups2(S4); # subgroups H of G up to conjugacy 
[ Group( () )^G, Group( [ (1,3)(2,4) ] )^G, Group( [ (3,4) ] )^G, 
  Group( [ (2,4,3) ] )^G, Group( [ (1,4)(2,3), (1,3)(2,4) ] )^G, 
  Group( [ (3,4), (1,2)(3,4) ] )^G, Group( [ (1,3,2,4), (1,2)(3,4) ] )^G, 
  Group( [ (3,4), (2,4,3) ] )^G, Group( [ (1,4)(2,3), (1,3)(2,4), (3,4) ] )^G,
  Group( [ (1,4)(2,3), (1,3)(2,4), (2,4,3) ] )^G, 
  Group( [ (1,4)(2,3), (1,3)(2,4), (2,4,3), (3,4) ] )^G ]
gap> Length(S4cs);
11
gap> S4H:=Hcandidates(S4); # exclude not core-free H=V4, D4, A4, S4
[ Group(()), Group([ (1,3)(2,4) ]), Group([ (3,4) ]), Group([ (2,4,3) ]), 
  Group([ (3,4), (1,2)(3,4) ]), Group([ (1,3,2,4), (1,2)(3,4) ]), 
  Group([ (3,4), (2,4,3) ]) ]
gap> Length(S4H);
7
gap> S4J:=List(S4H,x->Norm1TorusJCoset(S4,x));;
gap> List(S4J,IsInvertibleF); # T is not retract k-rational
[ false, false, false, false, false, false, false ]

gap> List(S4H,StructureDescription); # for checking core-free H
[ "1", "C2", "C2", "C3", "C2 x C2", "C4", "S3" ]

gap> S5:=SymmetricGroup(5); # G=S5=PGL(2,5) with |G|=120
Sym( [ 1 .. 5 ] )
gap> S5cs:=ConjugacyClassesSubgroups2(S5);; # subgroups H of G up to conjugacy 
gap> Length(S5cs);
19
gap> S5H:=Hcandidates(S5); # exclude not core-free H=A5, S5
[ Group(()), Group([ (1,3) ]), Group([ (2,3)(4,5) ]), Group([ (3,4,5) ]), 
  Group([ (2,3)(4,5), (2,4)(3,5) ]), Group([ (4,5), (2,3)(4,5) ]), 
  Group([ (2,5,3,4) ]), Group([ (1,2,3,4,5) ]), Group([ (3,4,5), (1,2) ]), 
  Group([ (3,4,5), (4,5) ]), Group([ (1,2)(4,5), (3,4,5) ]), 
  Group([ (2,3), (2,4,3,5) ]), Group([ (1,4)(2,3), (1,3)(4,5) ]), 
  Group([ (3,4,5), (2,4)(3,5) ]), Group([ (3,4,5), (4,5), (1,2)(4,5) ]), 
  Group([ (2,3,5,4), (1,4)(2,3) ]), Group([ (3,5,4), (2,5,4,3) ]) ]
gap> Length(S5H);
17
gap> S5J:=List(S5H,x->Norm1TorusJCoset(S5,x));;
gap> Filtered([1..17],x->IsInvertibleF(S5J[x])=true); 
# T is not retract k-rational except for i=5,12,14,17
[ 5, 12, 14, 17 ]
gap> List([5,12,14,17],x->S5H[x]);
[ Group([ (2,3)(4,5), (2,4)(3,5) ]), Group([ (2,3), (2,4,3,5) ]),
  Group([ (3,4,5), (2,4)(3,5) ]), Group([ (3,5,4), (2,5,4,3) ]) ]
gap> Filtered([1..Length(S5H)],x->StructureDescription(S5H[x])="C2 x C2");
[ 5, 6 ]
gap> StructureDescription(Intersection(DerivedSubgroup(S5),S5H[5]));
"C2 x C2"
gap> StructureDescription(Intersection(DerivedSubgroup(S5),S5H[6]));
"C2"
gap> IdSmallGroup(SymmetricGroup(5));
[ 120, 34 ]
gap> Filtered([1..NrTransitiveGroups(30)],x->Order(TransitiveGroup(30,x))=120);
[ 17, 18, 19, 20, 21, 22, 23, 24, 25, 26, 27, 28, 29, 30, 31, 32, 33, 34 ]
gap> Filtered(last,x->IdSmallGroup(TransitiveGroup(30,x))=[120,34]);
[ 22, 25, 27 ]
gap> List([22,25,27],x->StructureDescription(Intersection(Stabilizer(
> TransitiveGroup(30,x),1),DerivedSubgroup(TransitiveGroup(30,x)))));
[ "C2", "C2 x C2", "C2" ]
gap> StructureDescription(Intersection(DerivedSubgroup(S5),S5H[5])); 
# G=30T25 and H5=V4 with [G:H5]=30
"C2 x C2"

gap> F:=FlabbyResolution(S5J[5]).actionF; # [F]=[J_{G/H5}]^{fl} with H5=V4
<matrix group with 2 generators>
gap> Rank(F.1); # [F]=[J_{G/H5}]^{fl} with rank F=151
151
gap> ll:=PossibilityOfStablyPermutationF(S5J[5]);
[ [ 1, 0, 0, 0, 0, -2, 6, 1, -3, 3, -2, 4, 2, 5, 2, -6, -8, -3, 6, -2 ],
  [ 0, 1, 0, 0, 0, -1, -1, 0, 0, -1, 0, 0, 0, 0, 1, 1, 1, 0, -1, 0 ],
  [ 0, 0, 1, 0, 0, 0, -2, 0, 1, -1, 0, 0, -1, -1, 0, 2, 2, 1, -2, 0 ],
  [ 0, 0, 0, 1, 0, -2, 10, 1, -5, 5, -3, 4, 3, 6, 2, -10, -12, -4, 10, -2 ],
  [ 0, 0, 0, 0, 1, 2, -2, 0, 2, -2, 1, -2, -1, -2, -2, 2, 4, 1, -2, 0 ] ]
gap> List(ll,x->x[Length(x)]); # [F]\neq 0, i.e. T is not stably k-rational 
[ -2, 0, 0, -2, 0 ]

gap> F:=FlabbyResolution(S5J[12]).actionF; # [F]=[J_{G/H12}]^{fl} with H12=D4
<matrix group with 2 generators>
gap> Rank(F.1); # [F]=[J_{G/H12}]^{fl} with rank F=76
76
gap> ll:=PossibilityOfStablyPermutationF(S5J[12]);
[ [ 1, 0, 0, 0, 0, 0, -8, -1, 3, -5, 2, -2, -2, -5, 0, 8, 10, 3, -8, 2 ],
  [ 0, 1, 0, 0, 0, -1, -1, 0, 0, -1, 0, 0, 0, 0, 1, 1, 1, 0, -1, 0 ],
  [ 0, 0, 1, 0, 0, 0, -2, 0, 1, -1, 0, 0, -1, -1, 0, 2, 2, 1, -2, 0 ],
  [ 0, 0, 0, 1, 0, 0, -4, -1, 1, -3, 1, -2, -1, -4, 0, 4, 6, 2, -4, 2 ],
  [ 0, 0, 0, 0, 1, 2, -2, 0, 2, -2, 1, -2, -1, -2, -2, 2, 4, 1, -2, 0 ] ]
gap> List(ll,x->x[Length(x)]); # [F]\neq 0, i.e. T is not stably k-rational 
[ 2, 0, 0, 2, 0 ]

gap> F:=FlabbyResolution(S5J[14]).actionF; # [F]=[J_{G/H14}]^{fl} with H14=A4
<matrix group with 2 generators>
gap> Rank(F.1); # [F]=[J_{G/H14}]^{fl} with rank F=41
41
gap> ll:=PossibilityOfStablyPermutationF(S5J[14]);
[ [ 1, 0, 0, 0, 0, 0, -4, 1, 3, -3, 2, 0, -2, -1, 0, 4, 6, 1, -4, -2 ],
  [ 0, 1, 0, 0, 0, -1, -1, 0, 0, -1, 0, 0, 0, 0, 1, 1, 1, 0, -1, 0 ],
  [ 0, 0, 1, 0, 0, 0, -2, 0, 1, -1, 0, 0, -1, -1, 0, 2, 2, 1, -2, 0 ],
  [ 0, 0, 0, 1, 0, 0, 0, 1, 1, -1, 1, 0, -1, 0, 0, 0, 2, 0, 0, -2 ],
  [ 0, 0, 0, 0, 1, 2, -2, 0, 2, -2, 1, -2, -1, -2, -2, 2, 4, 1, -2, 0 ] ]
gap> List(ll,x->x[Length(x)]); # [F]\neq 0, i.e. T is not stably k-rational 
[ -2, 0, 0, -2, 0 ]

gap> F:=FlabbyResolution(S5J[17]).actionF; # [F]=[J_{G/H17}]^{fl} with H17=S4
<matrix group with 2 generators>
gap> Rank(F.1); # [F]=[J_{G/H17}]^{fl} with rank F=16
16
gap> ll:=PossibilityOfStablyPermutationF(S5J[17]);
[ [ 1, 0, 0, 0, 0, 0, -4, -1, 1, 0, -3, 0, 0, -1, 0, 4, 4, 1, -4, 2 ],
  [ 0, 1, 0, 0, 0, -1, -1, 0, 0, 0, -1, 0, 0, 0, 1, 1, 1, 0, -1, 0 ],
  [ 0, 0, 1, 0, 0, 0, -2, 0, 1, 0, -1, 0, -1, -1, 0, 2, 2, 1, -2, 0 ],
  [ 0, 0, 0, 1, 0, 0, 0, -1, -1, -1, -1, 0, 1, 0, 0, 0, 0, 0, 0, 2 ],
  [ 0, 0, 0, 0, 1, 2, -2, 0, 2, 1, -2, -2, -1, -2, -2, 2, 4, 1, -2, 0 ] ]
gap> List(ll,x->x[Length(x)]); # [F]\neq 0, i.e. T is not stably k-rational 
[ 2, 0, 0, 2, 0 ]

gap> List(S5H,StructureDescription); # for checking core-free H
[ "1", "C2", "C2", "C3", "C2 x C2", "C2 x C2", "C4", "C5", "C6", "S3", "S3",
  "D8", "D10", "A4", "D12", "C5 : C4", "S4" ]
gap> List([5,12,14,17],x->StructureDescription(S5H[x]));
[ "C2 x C2", "D8", "A4", "S4" ]

gap> S6:=SymmetricGroup(6); # G=S6 with |G|=720
Sym( [ 1 .. 6 ] )
gap> S6cs:=ConjugacyClassesSubgroups2(S6);; # subgroups H of G up to conjugacy 
gap> Length(S6cs);
56
gap> S6H:=Hcandidates(S6);; # exclude not core-free H=A6, S6
gap> Length(S6H);
54 
gap> S6J:=List(S6H,x->Norm1TorusJCoset(S6,x));;
gap> Filtered([1..54],x->IsInvertibleF(SylowSubgroup(S6J[x],2))=true);
[ 27, 34, 43, 44, 48, 49 ]
gap> List([27,34,43,44,48,49],x->IsInvertibleF(SylowSubgroup(S6J[x],3))); 
# T is not retract k-rational
[ false, false, false, false, false, false ]

gap> List(S6H,StructureDescription); # for checking core-free H
[ "1", "C2", "C2", "C2", "C3", "C3", "C2 x C2", "C2 x C2", "C2 x C2",
  "C2 x C2", "C2 x C2", "C4", "C4", "C5", "S3", "S3", "C6", "C6", "S3",
  "S3", "C2 x C2 x C2", "C2 x C2 x C2", "C4 x C2", "D8", "D8", "D8", "D8",
  "C3 x C3", "D10", "A4", "A4", "D12", "D12", "C2 x D8", "(C3 x C3) : C2",
  "C3 x S3", "C3 x S3", "C5 : C4", "C2 x A4", "S4", "C2 x A4", "S4", "S4",
  "S4", "S3 x S3", "S3 x S3", "(C3 x C3) : C4", "C2 x S4", "C2 x S4", "A5",
  "A5", "(S3 x S3) : C2", "S5", "S5" ]
\end{verbatim}~\\

\begin{verbatim}
gap> A4:=AlternatingGroup(4); # G=A4=PSL(2,3) with |G|=12
Alt( [ 1 .. 4 ] )
gap> A4cs:=ConjugacyClassesSubgroups2(A4); # subgroups H of G up to conjugacy 
[ Group( () )^G, Group( [ (1,2)(3,4) ] )^G, Group( [ (2,4,3) ] )^G, 
  Group( [ (1,3)(2,4), (1,2)(3,4) ] )^G, 
  Group( [ (1,3)(2,4), (1,2)(3,4), (2,4,3) ] )^G ]
gap> Length(A4cs);
5
gap> A4H:=Hcandidates(A4); # exclude not core-free H=V4, A4
[ Group(()), Group([ (1,2)(3,4) ]), Group([ (2,4,3) ]) ]
gap> Length(A4H);
3
gap> A4J:=List(A4H,x->Norm1TorusJCoset(A4,x));;
gap> List(A4J,IsInvertibleF); # T is not retract k-rational
[ false, false, false ]

gap> List(A4cs,x->StructureDescription(Representative(x))); # for checking H 
[ "1", "C2", "C3", "C2 x C2", "A4" ]
gap> List(A4H,StructureDescription); # for checking core-free H 
[ "1", "C2", "C3" ]

gap> A5:=AlternatingGroup(5); # G=A5=PSL(2,5)=PSL(2,4)=PGL(2,4)=SL(2,4) with |G|=60
Alt( [ 1 .. 5 ] )
gap> A5cs:=ConjugacyClassesSubgroups2(A5); # subgroups H of G up to conjugacy 
[ Group( () )^G, Group( [ (2,3)(4,5) ] )^G, Group( [ (3,4,5) ] )^G, 
  Group( [ (2,3)(4,5), (2,4)(3,5) ] )^G, Group( [ (1,2,3,4,5) ] )^G, 
  Group( [ (1,2)(4,5), (3,4,5) ] )^G, Group( [ (1,4)(2,3), (1,3)(4,5) ] )^G, 
  Group( [ (3,4,5), (2,4)(3,5) ] )^G, Group( [ (2,4)(3,5), (1,2,5) ] )^G ]
gap> Length(A5cs);
9
gap> A5H:=Hcandidates(A5); # exclude not core-free H=A5
[ Group(()), Group([ (2,3)(4,5) ]), Group([ (3,4,5) ]), Group([ (2,3)(4,5), (2,4)(3,5) ]), 
  Group([ (1,2,3,4,5) ]), Group([ (1,2)(4,5), (3,4,5) ]), Group([ (1,4)(2,3), (1,3)(4,5) ]), 
  Group([ (3,4,5), (2,4)(3,5) ]) ]
gap>Length(A5H);
8
gap> A5J:=List(A5H,x->Norm1TorusJCoset(A5,x));;
gap> List(A5J,IsInvertibleF); 
[ false, false, false, true, false, false, false, true ]
gap> Filtered([1..8],x->IsInvertibleF(A5J[x])=true); 
# T is not retract k-rational except for i=4,8
[ 4, 8 ]

gap> mi:=SearchCoflabbyResolutionBase(TransposedMatrixGroup(A5J[4]),0);; # H=C2xC2
gap> F:=FlabbyResolutionFromBase(A5J[4],mi).actionF;;
gap> Rank(F.1); # [F]=[J_{G/H4}]^{fl} with rank F=21
21
gap> ll:=PossibilityOfStablyPermutationFFromBase(A5J[4],mi);
[ [ 1, -2, -1, 0, 0, 1, 1, 1, -1, 0 ], [ 0, 0, 0, 0, 1, 0, 0, 2, -1, -1 ] ]
gap> l:=ll[2]; # possibility for Z[G/H5]+2Z[G/H8]=Z+F 
# with rank 12+2*5=1+21=22 where H5=C5,H8=A4
[ 0, 0, 0, 0, 1, 0, 0, 2, -1, -1 ]
gap> bp:=StablyPermutationFCheckPFromBase(A5J[4],mi,Nlist(l),Plist(l));;
gap> Length(bp);
16
gap> Length(bp[1]);
22
gap> rs:=RandomSource(IsMersenneTwister);
<RandomSource in IsMersenneTwister>
gap> rr:=List([1..200000],x->List([1..16],y->Random(rs,[-1..2])));;
gap> Filtered(rr,x->Determinant(x*bp)^2=1);
[ [ 2, 0, 1, 2, 1, 1, 1, 0, -1, 0, -1, -1, 1, -1, 2, 1 ] ]
gap> nn:=last[1];
[ 2, 0, 1, 2, 1, 1, 1, 0, -1, 0, -1, -1, 1, -1, 2, 1 ]
gap> P:=nn*bp;
[ [ 2, 2, 2, 2, 2, 2, 2, 2, 2, 2, 2, 2, 0, 0, 0, 0, 0, 1, 1, 1, 1, 1 ],
  [ 2, 1, 1, 1, 1, 1, 1, 1, 2, 2, 1, 1, 0, 0, 0, -1, 0, 0, 0, 0, -1, 0 ],
  [ 1, 2, 1, 1, 1, 2, 1, 1, 1, 1, 2, 1, 0, 0, -1, 0, 0, 0, 0, -1, 0, 0 ],
  [ 1, 1, 1, 2, 1, 1, 2, 1, 1, 1, 1, 2, -1, 0, 0, 0, 0, -1, 0, 0, 0, 0 ],
  [ 1, 1, 2, 1, 2, 1, 1, 2, 1, 1, 1, 1, 0, 0, 0, 0, -1, 0, 0, 0, 0, -1 ],
  [ -1, -1, -2, -1, -1, -1, -1, -1, -1, -2, -2, -2, 1, -1, 1, 1, 1, 2, 1, 2, 2, 2 ],
  [ 1, 1, 2, 1, 1, 2, 2, 1, 1, 1, 1, 1, 0, -1, 0, 0, 0, 0, -1, 0, 0, 0 ],
  [ -2, -1, -1, -1, -1, -1, -2, -1, -1, -1, -2, -2, 1, 1, 1, -1, 1, 2, 2, 2, 1, 2 ],
  [ -1, -2, -1, -1, -2, -2, -1, -1, -2, -1, -1, -1, 1, 1, 1, -1, 1, 2, 2, 2, 1, 2 ],
  [ -2, -1, -2, -1, -1, -1, -2, -1, -1, -2, -1, -1, -1, 1, 1, 1, 1, 1, 2, 2, 2, 2 ],
  [ -2, -1, -1, -1, -1, -2, -2, -1, -2, -1, -1, -1, 1, 1, -1, 1, 1, 2, 2, 1, 2, 2 ],
  [ -1, -1, -1, -2, -1, -1, -1, -2, -1, -1, -2, -2, 1, 1, -1, 1, 1, 2, 2, 1, 2, 2 ],
  [ -1, -1, -2, -1, -1, -2, -1, -1, -2, -2, -1, -1, 1, 1, 1, 1, -1, 2, 2, 2, 2, 1 ],
  [ -1, -2, -1, -1, -2, -1, -1, -1, -1, -1, -2, -2, 1, 1, 1, 1, -1, 2, 2, 2, 2, 1 ],
  [ -2, -1, -1, -2, -1, -1, -2, -2, -1, -1, -1, -1, 1, 1, 1, 1, -1, 2, 2, 2, 2, 1 ],
  [ 4, 4, 3, 4, 4, 4, 4, 4, 4, 3, 5, 5, -3, -3, -1, -1, -1, -6, -6, -5, -5, -5 ], 
  [ 1, 1, 1, 1, 2, 1, 1, 1, 2, 1, 1, 2, 0, -1, 0, 0, 0, 0, -1, 0, 0, 0 ],
  [ 2, 1, 1, 1, 1, 1, 1, 2, 1, 1, 2, 1, 0, -1, 0, 0, 0, 0, -1, 0, 0, 0 ],
  [ -1, -1, -1, -1, -2, -1, -1, -2, -1, 0, -2, -2, 0, 1, 0, 0, 0, 0, 1, 0, 0, 0 ],
  [ 1, 2, 1, 1, 1, 1, 1, 1, 1, 2, 2, 1, -1, 0, 0, 0, 0, -1, 0, 0, 0, 0 ],
  [ 2, 1, 1, 1, 1, 1, 1, 1, 1, 2, 2, 1, 0, 0, 0, 0, -1, 0, 0, 0, 0, -1 ],
  [ 1, 1, 1, 2, 1, 1, 1, 1, 2, 1, 1, 2, 0, 0, 0, 0, -1, 0, 0, 0, 0, -1 ] ]
gap> Determinant(P);
-1
gap> StablyPermutationFCheckMatFromBase(A5J[4],mi,Nlist(l),Plist(l),P);
# T is stably k-rational
true

gap> F:=FlabbyResolutionLowRank(A5J[8]).actionF; # H=A4
<matrix group with 2 generators>
gap> Rank(F.1); # [F]=[J_{G/H8}]^{fl} with rank F=16
16
gap> mi:=SearchCoflabbyResolutionBaseLowRank(TransposedMatrixGroup(A5J[8]),0);;
gap> ll:=PossibilityOfStablyPermutationFFromBase(A5J[8],mi);
[ [ 1, -2, -1, 0, 0, 1, 1, 1, -1, 0 ], [ 0, 0, 0, 0, 1, 1, -1, 0, 0, -1 ] ]
gap> l:=ll[2]; # possibility for Z[G/H5]+Z[G/H6]=Z[G/H7]+F 
# with rank 12+10=6+16=22 where H5=C5,H6=S3,H7=D5
[ 0, 0, 0, 0, 1, 1, -1, 0, 0, -1 ]
gap> bp:=StablyPermutationFCheckPFromBase(A5J[8],mi,Nlist(l),Plist(l));; 
gap> Length(bp);
11
gap> Length(bp[1]);
22
gap> SearchPRowBlocks(bp);
rec( bpBlocks := [ [ 1, 2, 3, 4 ], [ 5, 6, 7, 8, 9, 10, 11 ] ], 
  rowBlocks := 
    [ [ 1, 2, 3, 4, 5, 6 ], 
      [ 7, 8, 9, 10, 11, 12, 13, 14, 15, 16, 17, 18, 19, 20, 21, 22 ] ] )
gap> r1:=SearchPFilterRowBlocks(bp,[1..4],[1..6],3,[-1..2]);;
gap> Length(r1);
94
gap> r2:=SearchPFilterRowBlocks(bp,[5..11],[7..22],2);;
gap> Length(r2);
8
gap> P:=SearchPMergeRowBlock(r1,r2);
[ [ 2, 1, 1, 1, 1, 2, 1, 1, 1, 1, 1, 1, 1, 1, 0, 0, 1, 0, 0, 1, 1, 0 ], 
  [ 1, 2, 1, 1, 1, 1, 1, 1, 1, 1, 2, 1, 0, 1, 1, 0, 0, 1, 1, 0, 1, 0 ], 
  [ 1, 1, 2, 1, 1, 1, 2, 1, 1, 1, 1, 1, 1, 0, 1, 1, 0, 0, 0, 0, 1, 1 ], 
  [ 1, 1, 1, 2, 1, 1, 1, 1, 1, 2, 1, 1, 0, 0, 1, 1, 1, 0, 1, 1, 0, 0 ], 
  [ 1, 1, 1, 1, 2, 1, 1, 1, 2, 1, 1, 1, 0, 1, 0, 1, 0, 1, 0, 1, 0, 1 ], 
  [ 1, 1, 1, 1, 1, 1, 1, 2, 1, 1, 1, 2, 1, 0, 0, 0, 1, 1, 1, 0, 0, 1 ], 
  [ 0, 0, 1, 1, 1, 0, 0, 0, 0, 0, 0, 0, 0, 0, 0, 1, 0, 0, 0, 0, 0, 0 ], 
  [ 0, 0, 0, 0, 0, 1, 1, 1, 0, 0, 0, 0, 1, 0, 0, 0, 0, 0, 0, 0, 0, 0 ], 
  [ 1, 1, 0, 0, 0, 0, 0, 0, 1, 0, 0, 0, 0, 1, 0, 0, 0, 0, 0, 0, 0, 0 ], 
  [ 0, 0, 0, 0, 0, 0, 0, 0, 0, 1, 1, 1, 0, 0, 0, 0, 0, 0, 1, 0, 0, 0 ], 
  [ 0, 0, 0, 0, 0, 1, 1, 0, 0, 0, 1, 0, 0, 0, 0, 0, 0, 0, 0, 0, 1, 0 ], 
  [ 1, 0, 0, 0, 0, 0, 0, 0, 1, 1, 0, 0, 0, 0, 0, 0, 0, 0, 0, 1, 0, 0 ], 
  [ 0, 0, 1, 0, 1, 0, 0, 0, 0, 0, 0, 1, 0, 0, 0, 0, 0, 0, 0, 0, 0, 1 ], 
  [ 1, 0, 0, 0, 0, 0, 0, 0, 0, 1, 0, 1, 0, 0, 0, 0, 1, 0, 0, 0, 0, 0 ], 
  [ 0, 1, 1, 1, 0, 0, 0, 0, 0, 0, 0, 0, 0, 0, 1, 0, 0, 0, 0, 0, 0, 0 ], 
  [ 0, 1, 0, 0, 0, 0, 0, 1, 1, 0, 0, 0, 0, 0, 0, 0, 0, 1, 0, 0, 0, 0 ], 
  [ 0, 0, 0, 0, 1, 0, 0, 0, 0, 0, 1, 1, 0, 0, 0, 0, 0, 1, 0, 0, 0, 0 ], 
  [ 0, 0, 0, 0, 0, 0, 1, 0, 0, 1, 1, 0, 0, 0, 1, 0, 0, 0, 0, 0, 0, 0 ], 
  [ 0, 0, 0, 1, 0, 1, 0, 1, 0, 0, 0, 0, 0, 0, 0, 0, 1, 0, 0, 0, 0, 0 ], 
  [ 0, 0, 0, 0, 0, 0, 1, 1, 1, 0, 0, 0, 0, 0, 0, 0, 0, 0, 0, 0, 0, 1 ], 
  [ 0, 0, 0, 1, 1, 1, 0, 0, 0, 0, 0, 0, 0, 0, 0, 0, 0, 0, 0, 1, 0, 0 ], 
  [ 1, 1, 1, 0, 0, 0, 0, 0, 0, 0, 0, 0, 0, 0, 0, 0, 0, 0, 0, 0, 1, 0 ] ]
gap> Determinant(P);
1
gap> StablyPermutationFCheckMatFromBase(A5J[8],mi,Nlist(l),Plist(l),P);
# T is stably k-rational
true

gap> List(A5H,StructureDescription); # for checking core-free H
[ "1", "C2", "C3", "C2 x C2", "C5", "S3", "D10", "A4" ]
gap> List([4,8],x->StructureDescription(A5H[x])); 
[ "C2 x C2", "A4" ]

gap> A6:=AlternatingGroup(6); # G=A6=PSL(2,9) with |G|=360
Alt( [ 1 .. 6 ] )
gap> A6cs:=ConjugacyClassesSubgroups2(A6);; # subgroups H of G up to conjugacy 
gap> Length(A6cs); 
22
gap> A6H:=Hcandidates(A6); # exclude not core-free H=A6
gap> Length(A6H); 
21
gap> A6J:=List(A6H,x->Norm1TorusJCoset(A6,x));;
gap> Filtered([1..21],x->IsInvertibleF(SylowSubgroup(A6J[x],2))=true);
[ 11, 17, 18 ]
gap> List([11,17,18],x->IsInvertibleF(SylowSubgroup(A6J[x],3))); 
# T is not retract k-rational
[ false, false, false ]

gap> List(A6H,StructureDescription); # for checking core-free H
[ "1", "C2", "C3", "C3", "C2 x C2", "C2 x C2", "C4", "C5", "S3", "S3", "D8",
  "C3 x C3", "D10", "A4", "A4", "(C3 x C3) : C2", "S4", "S4",
  "(C3 x C3) : C4", "A5", "A5" ]
\end{verbatim}~\\

\begin{verbatim}
gap> GL23:=GL(2,3); # G=GL(2,3) with |G|=48 
GL(2,3)
gap> GL23cs:=ConjugacyClassesSubgroups2(GL23);; # subgroups H of G up to conjugacy 
gap> Length(GL23cs);
16 
gap> GL23H:=Hcandidates(GL23);; # exclude not core-free H with Z(GL(2,3))=C2<=H 
gap> Length(GL23H); 
5
gap> GL23J:=List(GL23H,x->Norm1TorusJCoset(GL23,x));;
gap> List(GL23J,IsInvertibleF); # T is not retract k-rational
[ false, false, false, false, false ]

gap> List(GL23cs,x->StructureDescription(Representative(x))); # for checking H
[ "1", "C2", "C2", "C3", "C4", "C2 x C2", "S3", "S3", "C6", "Q8", "D8",
  "C8", "D12", "QD16", "SL(2,3)", "GL(2,3)" ]
gap> List(GL23H,StructureDescription); # for checking core-free H
[ "1", "C2", "C3", "S3", "S3" ]

gap> GL24:=GL(2,4); # G=GL(2,4)=A5xC3 with |G|=180
GL(2,4)
gap> GL24cs:=ConjugacyClassesSubgroups2(GL24);; # subgroups H of G up to conjugacy 
gap> Length(GL24cs);
21
gap> GL24H:=Hcandidates(GL24);; # exclude not core-free H with Z(GL(2,4))=C3<=H and H=SL(2,4)
gap> Length(GL24H);
11
gap> GL24J:=List(GL24H,x->Norm1TorusJCoset(GL24,x));;
gap> Filtered([1..11],x->IsInvertibleF(SylowSubgroup(GL24J[x],2))=true); 
# T is not retract k-rational except for i=5,11
[ 5, 11 ]
gap> IsInvertibleF(SylowSubgroup(GL24J[5],3)); # T is not retract k-rational
false
gap> IsInvertibleF(SylowSubgroup(GL24J[11],3));
true
gap> IsInvertibleF(GL24J[11]);
true
gap> Filtered([1..NrTransitiveGroups(15)],x->Order(TransitiveGroup(15,x))=180);
[ 15, 16 ]
gap> G15:=TransitiveGroup(15,15); # G15=15T15
3A_5(15)=[3]A(5)=GL(2,4)
gap> G16:=TransitiveGroup(15,16); # G16=15T16
A(5)[x]3
gap> IdSmallGroup(G15); # G15=A5xC3
[ 180, 19 ]
gap> IdSmallGroup(G16); # G16=A5xC3
[ 180, 19 ]
gap> IdSmallGroup(GL24); # GL24=A5xC3
[ 180, 19 ]
gap> StructureDescription(Intersection(DerivedSubgroup(G15),Stabilizer(G15,1)));
"C2 x C2"
gap> StructureDescription(Intersection(DerivedSubgroup(G16),Stabilizer(G16,1)));
"A4"
gap> StructureDescription(Intersection(DerivedSubgroup(GL24),GL24H[9])); # G/H9=15T15
"C2 x C2"
gap> StructureDescription(Intersection(DerivedSubgroup(GL24),GL24H[10])); # G/H10=15T15
"C2 x C2"
gap> StructureDescription(Intersection(DerivedSubgroup(GL24),GL24H[11])); # G/H11=15T16
"A4"
gap> StructureDescription(DerivedSubgroup(GL24));
"A5"
gap> F:=FlabbyResolutionLowRank(Norm1TorusJTransitiveGroup(15,16)).actionF;
<matrix group with 2 generators>
gap> Rank(F.1); # [F]=[J_{G/H11}]^{fl} with rank F=36
36
gap> F2:=FlabbyResolutionLowRankFromGroup(F,Norm1TorusJTransitiveGroup(15,16)).actionF;
<matrix group with 2 generators>
gap> Rank(F2.1);
14
gap> F3:=FlabbyResolutionLowRankFromGroup(F2,Norm1TorusJTransitiveGroup(15,16)).actionF; 
# T is stably k-rational 
Group([ [ [ 1 ] ], [ [ 1 ] ] ])

gap> List(GL24cs,x->StructureDescription(Representative(x))); # for checking H
[ "1", "C2", "C3", "C3", "C3", "C2 x C2", "C5", "S3", "C6", "C3 x C3",
  "D10", "C6 x C2", "A4", "A4", "A4", "C15", "C3 x S3", "C3 x D10",
  "C3 x A4", "A5", "GL(2,4)" ]
gap> List(GL24H,StructureDescription); # for checking core-free H
[ "1", "C2", "C3", "C3", "C2 x C2", "C5", "S3", "D10", "A4", "A4", "A4" ]

gap> GL25:=GL(2,5);
GL(2,5)
gap> GL25cs:=ConjugacyClassesSubgroups2(GL25);; # subgroups H of G up to conjugacy 
gap> Length(GL25cs);
48
gap> GL25H:=Hcandidates(GL25);; # exclude not core-free H with Z(GL(2,5))\cap H\neq 1 
gap> Length(GL25H);
13
gap> GL25J:=List(GL25H,x->Norm1TorusJCoset(GL25,x));;
gap> Filtered([1..13],x->IsInvertibleF(SylowSubgroup(GL25J[x],2))=true); 
# T is not retract k-rational 
[  ]

gap> List(GL25cs,x->StructureDescription(Representative(x))); # for checking H
[ "1", "C2", "C2", "C3", "C4", "C2 x C2", "C4", "C4", "C4", "C5", "C6",
  "S3", "Q8", "C8", "D8", "C4 x C2", "C4 x C2", "C10", "D10", "D10", "D12",
  "C3 : C4", "C12", "(C4 x C2) : C2", "C4 x C4", "C8 : C2", "C5 : C4",
  "C5 : C4", "D20", "C5 : C4", "C5 : C4", "C5 : C4", "C20", "SL(2,3)",
  "C4 x S3", "C3 : C8", "C24", "(C4 x C4) : C2", "C2 x (C5 : C4)",
  "C4 x D10", "C2 x (C5 : C4)", "((C4 x C2) : C2) : C3", "C24 : C2",
  "C4 x (C5 : C4)", "SL(2,3) : C4", "SL(2,5)", "SL(2,5) : C2", "GL(2,5)" ]
gap> List(GL25H,StructureDescription); # for checking core-free H
[ "1", "C2", "C3", "C4", "C4", "C5", "S3", "D10", "D10", "C5 : C4",
  "C5 : C4", "C5 : C4", "C5 : C4" ]
\end{verbatim}~\\

\begin{verbatim}
gap> SL23:=SL(2,3); # G=SL(2,3) with |G|=24 
SL(2,3)
gap> SL23cs:=ConjugacyClassesSubgroups2(SL23);; # subgroups H of G up to conjugacy 
gap> Length(SL23cs);
7
gap> SL23H:=Hcandidates(SL23); # exclude not core-free H with Z(SL(2,3))=C2<=H
[ Group([  ]), Group([ [ [ 0*Z(3), Z(3) ], [ Z(3)^0, Z(3) ] ] ]) ]
gap> Length(SL23H);
2
gap> SL23J:=List(SL23H,x->Norm1TorusJCoset(SL23,x));;
gap> List(SL23J,IsInvertibleF); # T is not retract k-rational
[ false, false ]

gap> List(SL23cs,x->StructureDescription(Representative(x))); # for checking H
[ "1", "C2", "C3", "C4", "C6", "Q8", "SL(2,3)" ]
gap> List(SL23H,StructureDescription); # for checking core-free H
[ "1", "C3" ]

gap> SL25:=SL(2,5); # G=SL(2,5) with |G|=120 
SL(2,5)
gap> SL25cs:=ConjugacyClassesSubgroups2(SL25);; # subgroups H of G up to conjugacy 
gap> Length(SL25cs); 
12
gap> SL25H:=Hcandidates(SL25); # exclude not core-free H with Z(SL(2,5))=C2<=H
[ Group([  ]), Group([ [ [ 0*Z(5), Z(5) ], [ Z(5), Z(5)^2 ] ] ]), 
  Group([ [ [ 0*Z(5), Z(5)^2 ], [ Z(5)^0, Z(5) ] ] ]) ]
gap> Length(SL25H);
3
gap> SL25J:=List(SL25H,x->Norm1TorusJCoset(SL25,x));;
gap> List(SL25J,IsInvertibleF); # T is not retract k-rational
[ false, false, false ]

gap> List(SL25cs,x->StructureDescription(Representative(x))); # for checking H
[ "1", "C2", "C3", "C4", "C5", "C6", "Q8", "C10", "C3 : C4", "C5 : C4",
  "SL(2,3)", "SL(2,5)" ]
gap> List(SL25H,StructureDescription); # for checking core-free H
[ "1", "C3", "C5" ]

gap> SL27:=SL(2,7);
SL(2,7)
gap> SL27cs:=ConjugacyClassesSubgroups2(SL27);;
gap> Length(SL27cs);
19
gap> SL27H:=Hcandidates(SL27);; # exclude not core-free H with Z(SL(2,7))=C2<=H
gap> Length(SL27H);
4
gap> SL27J:=List(SL27H,x->Norm1TorusJCoset(SL27,x));;
gap> Filtered([1..4],x->IsInvertibleF(SylowSubgroup(SL27J[x],2))=true);
# T is not retract k-rational 
[  ]

gap> List(SL27cs,x->StructureDescription(Representative(x))); # for checking H
[ "1", "C2", "C3", "C4", "C6", "C7", "Q8", "Q8", "C8", "C3 : C4", "C14",
  "Q16", "C7 : C3", "SL(2,3)", "SL(2,3)", "C2 x (C7 : C3)",
  "C2 . S4 = SL(2,3) . C2", "C2 . S4 = SL(2,3) . C2", "SL(2,7)" ]
gap> List(SL27H,StructureDescription); # for checking core-free H
[ "1", "C3", "C7", "C7 : C3" ]
\end{verbatim}~\\

\begin{verbatim}
gap> PSL27:=PSL(2,7); # G=PSL(2,7)=PSL(3,2) with |G|=168
Group([ (3,7,5)(4,8,6), (1,2,6)(3,4,8) ])
gap> PSL27cs:=ConjugacyClassesSubgroups2(PSL27);; # subgroups H of G up to conjugacy 
gap> Length(PSL27cs);
15 
gap> PSL27H:=Hcandidates(PSL27);; # exclude not core-free H=PSL(2,7)
gap> Length(PSL27H);
14 
gap> PSL27J:=List(PSL27H,x->Norm1TorusJCoset(PSL27,x));;
gap> Filtered([1..14],x->IsInvertibleF(PSL27J[x])=true); 
# T is not retract k-rational except for i=9,13,14
[ 9, 13, 14 ]
gap> Filtered([1..NrTransitiveGroups(21)],x->Order(TransitiveGroup(21,x))=168); 
# G=21T14 and H=H9 with [G:H]=21
[ 14 ]

gap> F:=FlabbyResolution(PSL27J[9]).actionF; # [F]=[J_{G/H9}]^{fl} with H9=D4
<matrix group with 2 generators>
gap> Rank(F.1); # [F]=[J_{G/H9}]^{fl} with rank F=148
148
gap> ll:=PossibilityOfStablyPermutationF(PSL27J[9]);
[ [ 1, 0, 1, 0, 1, 0, 4, 2, 8, 4, 3, -4, -6, -4, 2, -4 ],
  [ 0, 1, 1, 0, 1, 0, 1, 1, 3, 2, 1, -2, -3, -1, 1, -2 ],
  [ 0, 0, 2, 0, 1, 0, 2, 1, 4, 2, 1, -3, -4, -2, 2, -2 ],
  [ 0, 0, 0, 1, -1, 0, 0, 0, 0, -1, 1, 0, 2, -2, 0, 0 ] ]
gap> List(ll,x->x[Length(x)]); # [F]\neq 0, i.e. T is not stably k-rational 
[ -4, -2, -2, 0 ]

gap> F:=FlabbyResolution(PSL27J[13]).actionF; # [F]=[J_{G/H13}]^{fl} with H13=S4
<matrix group with 2 generators>
gap> Rank(F.1); # [F]=[J_{G/H13}]^{fl} with rank F=36
36
gap> ll:=PossibilityOfStablyPermutationF(PSL27J[13]);
[ [ 1, 0, -3, 0, 0, 0, 0, 1, 0, 0, 2, 1, 2, 0, -2, -2 ],
  [ 0, 1, -1, 0, 0, 0, -1, 0, -1, 0, 0, 1, 1, 1, -1, 0 ],
  [ 0, 0, 0, 1, 0, 0, 0, 1, 0, -1, 2, -1, 2, -2, 0, -2 ],
  [ 0, 0, 0, 0, 1, 0, 0, 1, 0, 0, 1, -1, 0, 0, 0, -2 ] ]
gap> List(ll,x->x[Length(x)]); # [F]\neq 0, i.e. T is not stably k-rational 
[ -2, 0, -2, -2 ]

gap> F:=FlabbyResolution(PSL27J[14]).actionF; # [F]=[J_{G/H14}]^{fl} with H14=S4
<matrix group with 2 generators>
gap> Rank(F.1); # [F]=[J_{G/H14}]^{fl} with rank F=36
36
gap> ll:=PossibilityOfStablyPermutationF(PSL27J[14]);
[ [ 1, 0, -3, 0, 0, 0, 0, 1, 0, 0, 2, 1, 2, 0, -2, -2 ],
  [ 0, 1, -1, 0, 0, 0, -1, 0, -1, 0, 0, 1, 1, 1, -1, 0 ],
  [ 0, 0, 0, 1, 0, 0, 0, 1, 0, -1, 2, -1, 2, -2, 0, -2 ],
  [ 0, 0, 0, 0, 1, 0, 0, 1, 0, 0, 1, -1, 0, 0, 0, -2 ] ]
gap> List(ll,x->x[Length(x)]); # [F]\neq 0, i.e. T is not stably k-rational 
[ -2, 0, -2, -2 ]

gap> List(PSL27H,StructureDescription); # for checking core-free H
[ "1", "C2", "C3", "C2 x C2", "C2 x C2", "C4", "S3", "C7", "D8", "A4", "A4",
  "C7 : C3", "S4", "S4" ]
gap> List([9,13,14],x->StructureDescription(PSL27H[x]));
[ "D8", "S4", "S4" ] 
\end{verbatim}~\\

\begin{verbatim}
gap> PSL28:=PSL(2,8); # G=PSL(2,8) with |G|=504
Group([ (3,8,6,4,9,7,5), (1,2,3)(4,7,5)(6,9,8) ])
gap> PSL28cs:=ConjugacyClassesSubgroups2(PSL28); # subgroups H of G up to conjugacy 
[ Group( () )^G, Group( [ (1,5)(2,3)(6,8)(7,9) ] )^G,
  Group( [ (1,6,5)(2,3,9)(4,7,8) ] )^G,
  Group( [ (1,9)(3,7)(4,5)(6,8), (1,6)(3,5)(4,7)(8,9) ] )^G,
  Group( [ (1,3,2)(4,5,7)(6,8,9), (1,2)(4,9)(5,8)(6,7) ] )^G,
  Group( [ (1,9,6,4,7,8,3) ] )^G,
  Group( [ (1,9)(3,7)(4,5)(6,8), (1,4)(3,8)(5,9)(6,7), (1,6)(3,5)(4,7)(8,9) ] )^G,
  Group( [ (1,5,2,7,8,4,6,9,3) ] )^G,
  Group( [ (2,4)(3,6)(5,7)(8,9), (1,6)(3,5)(4,7)(8,9) ] )^G,
  Group( [ (1,2)(4,9)(5,8)(6,7), (1,6)(2,8)(3,9)(4,5) ] )^G,
  Group( [ (1,4)(3,8)(5,9)(6,7), (1,6,7,3,9,4,8) ] )^G,
  Group( [ (1,5)(2,3)(6,8)(7,9), (1,6,5)(2,3,9)(4,7,8) ] )^G ]
gap> Length(PSL28cs); 
12 
gap> PSL28H:=Hcandidates(PSL28); # exclude not core-free H=PSL(2,8)
[ Group(()), Group([ (1,5)(2,3)(6,8)(7,9) ]), 
  Group([ (1,6,5)(2,3,9)(4,7,8) ]), 
  Group([ (1,9)(3,7)(4,5)(6,8), (1,6)(3,5)(4,7)(8,9) ]), 
  Group([ (1,3,2)(4,5,7)(6,8,9), (1,2)(4,9)(5,8)(6,7) ]), 
  Group([ (1,9,6,4,7,8,3) ]), 
  Group([ (1,9)(3,7)(4,5)(6,8), (1,4)(3,8)(5,9)(6,7), (1,6)(3,5)(4,7)(8,9) ]), 
  Group([ (1,5,2,7,8,4,6,9,3) ]), 
  Group([ (2,4)(3,6)(5,7)(8,9), (1,6)(3,5)(4,7)(8,9) ]), 
  Group([ (1,2)(4,9)(5,8)(6,7), (1,6)(2,8)(3,9)(4,5) ]), 
  Group([ (1,4)(3,8)(5,9)(6,7), (1,6,7,3,9,4,8) ]) ]
gap> Length(PSL28H);
11
gap> PSL28J:=List(PSL28H,x->Norm1TorusJCoset(PSL28,x));;
gap> Filtered([1..11],x->IsInvertibleF(SylowSubgroup(PSL28J[x],2))=true); 
# T is not retract k-rational except for i=7,11
[ 7, 11 ]

gap> mi:=SearchCoflabbyResolutionBaseLowRank(TransposedMatrixGroup(PSL28J[7]),0);; 
# H=C2xC2xC2
gap> F:=FlabbyResolutionFromBase(PSL28J[7],mi).actionF;;
gap> Rank(F.1); # [F]=[J_{G/H7}]^{fl} with rank F=73
73
gap> ll:=PossibilityOfStablyPermutationFFromBase(PSL28J[7],mi);
[ [ 1, -2, 0, 0, 0, -1, 0, 0, 1, 1, 1, -1, 0 ], 
  [ 0, 0, 0, 0, 0, 0, 0, 1, 0, 0, 2, -1, -1 ] ]
gap> l:=ll[2]; 
# possibility for Z[G/H8]+2Z[G/H11]=Z+F 
# with rank 56+2x9=1+73=74 where H8=C9,H11=(C2xC2xC2):C7
[ 0, 0, 0, 0, 0, 0, 0, 1, 0, 0, 2, -1, -1 ]
gap> I13:=IdentityMat(13);;
gap> bpS:=StablyPermutationFCheckPFromBasePari(PSL28J[7],mi,
> Nlist(l)+I13[5],Plist(l)+I13[5]:parisize:=2^33);;
# we need at least 8GB memory for PARI/GP computations
# possibility for Z[G/H5]+Z[G/H8]+2Z[G/H11]=Z[G/H5]+Z+F
# with rank 84+56+2x9=84+1+73=158 where H5=S3,H8=C9,H11=(C2xC2xC2):C7
# bpS is available from the web as PSL28E8bpS.dat 
gap> SearchPRowBlocks(bpS);
rec( 
  bpBlocks := 
    [ [ 1, 2, 3, 4, 5, 6, 7, 8, 9, 10, 11, 12, 13, 14, 15, 16, 17, 18, 19, 
          20, 21, 22, 23, 24, 25, 26, 27, 28, 29, 30, 31 ], 
      [ 32, 33, 34, 35 ], 
      [ 36, 37, 38, 39, 40, 41, 42, 43, 44, 45, 46, 47, 48, 49, 50, 51, 52, 
          53, 54, 55, 56, 57, 58, 59, 60, 61, 62, 63, 64, 65 ] ], 
  rowBlocks := 
    [ [ 1, 2, 3, 4, 5, 6, 7, 8, 9, 10, 11, 12, 13, 14, 15, 16, 17, 18, 19, 
          20, 21, 22, 23, 24, 25, 26, 27, 28, 29, 30, 31, 32, 33, 34, 35, 36, 
          37, 38, 39, 40, 41, 42, 43, 44, 45, 46, 47, 48, 49, 50, 51, 52, 53, 
          54, 55, 56, 57, 58, 59, 60, 61, 62, 63, 64, 65, 66, 67, 68, 69, 70, 
          71, 72, 73, 74, 75, 76, 77, 78, 79, 80, 81, 82, 83, 84 ], [ 85 ], 
      [ 86, 87, 88, 89, 90, 91, 92, 93, 94, 95, 96, 97, 98, 99, 100, 101, 
          102, 103, 104, 105, 106, 107, 108, 109, 110, 111, 112, 113, 114, 
          115, 116, 117, 118, 119, 120, 121, 122, 123, 124, 125, 126, 127, 
          128, 129, 130, 131, 132, 133, 134, 135, 136, 137, 138, 139, 140, 
          141, 142, 143, 144, 145, 146, 147, 148, 149, 150, 151, 152, 153, 
          154, 155, 156, 157, 158 ] ] )
gap> r1:=SearchPFilterRowBlocks(bpS,[1..31],[1..84],4);;
gap> Length(r1);
3100
gap> r2:=SearchPFilterRowBlocks(bpS,[32..35],[85],4,[-8..9]);;
gap> Length(r2);
97040
gap> r3:=SearchPFilterRowBlocks(bpS,[36..65],[86..158],4);;
gap> Length(r3);
264
gap> r13:=SearchPMergeRowBlock(r1,r3);;
gap> Length(r13);
8
gap> r132:=SearchPMergeRowBlock(r13,r2);;
gap> P1:=r132{[1..84]};;
gap> P2:=r132{[158]};;
gap> P3:=r132{[85..157]};;
gap> P:=Concatenation([P1,P2,P3]);;
gap> Size(P); # 158x158 matrix
158
gap> Determinant(P);
-1
gap> StablyPermutationFCheckMatFromBase(PSL28J[7],mi,Nlist(l)+I13[5],Plist(l)+I13[5],P);
# T is stably k-rational
true

gap> mi:=SearchCoflabbyResolutionBaseLowRank(TransposedMatrixGroup(PSL28J[11]),0);;
# H=(C2xC2xC2):C7
gap> F:=FlabbyResolutionFromBase(PSL28J[11],mi).actionF;;
gap> Rank(F.1); # [F]=[J_{G/H11}]^{fl} with rank F=64
64
gap> ll:=PossibilityOfStablyPermutationFFromBase(PSL28J[11],mi);
[ [ 1, -2, 0, 0, 0, -1, 0, 0, 1, 1, 1, -1, 0 ], 
  [ 0, 0, 0, 0, 0, 0, 0, 1, 1, -1, 0, 0, -1 ] ]
gap> l:=ll[2]; # possibility for Z[G/H8]+Z[G/H9]=Z[G/H10]+F 
# with rank 56+36=28+64=92 where H8=C9,H9=D7,H10=D9
[ 0, 0, 0, 0, 0, 0, 0, 1, 1, -1, 0, 0, -1 ]
gap> I13:=IdentityMat(13);;
gap> bpQ:=StablyPermutationFCheckPFromBasePari(PSL28J[11],mi,
> Nlist(l)+I13[3],Plist(l)+I13[3]:parisize:=2^34);; 
# we need at least 16GB memory for PARI/GP computations
# possibility for Z[G/H3]+Z[G/H8]+Z[G/H9]=Z[G/H3]+Z[G/H10]+F 
# with rank 168+56+36=168+28+64=260 where H3=C3,H8=C9,H9=D7,H10=D9
# bpQ is available from the web as PSL28bpQ.dat 
gap> SearchPRowBlocks(bpQ);
rec( 
  bpBlocks := 
    [ [ 1, 2, 3, 4, 5, 6, 7, 8, 9, 10, 11, 12, 13, 14, 15, 16, 17, 18, 19, 
          20, 21, 22, 23, 24, 25, 26, 27, 28, 29, 30, 31, 32, 33, 34, 35, 36, 
          37, 38, 39, 40, 41, 42, 43, 44, 45, 46, 47, 48, 49, 50, 51, 52, 53, 
          54, 55, 56, 57, 58, 59, 60, 61, 62, 63, 64, 65, 66, 67, 68, 69, 70, 
          71, 72, 73, 74, 75, 76, 77, 78, 79, 80, 81, 82, 83, 84, 85, 86, 87, 
          88, 89, 90, 91, 92 ], 
      [ 93, 94, 95, 96, 97, 98, 99, 100, 101, 102, 103, 104, 105, 106, 107, 
          108, 109, 110 ], 
      [ 111, 112, 113, 114, 115, 116, 117, 118, 119, 120, 121, 122, 123, 124, 
          125, 126, 127, 128, 129, 130, 131, 132, 133, 134, 135, 136, 137, 
          138, 139, 140, 141, 142, 143, 144, 145 ] ], 
  rowBlocks := 
    [ [ 1, 2, 3, 4, 5, 6, 7, 8, 9, 10, 11, 12, 13, 14, 15, 16, 17, 18, 19, 
          20, 21, 22, 23, 24, 25, 26, 27, 28, 29, 30, 31, 32, 33, 34, 35, 36, 
          37, 38, 39, 40, 41, 42, 43, 44, 45, 46, 47, 48, 49, 50, 51, 52, 53, 
          54, 55, 56, 57, 58, 59, 60, 61, 62, 63, 64, 65, 66, 67, 68, 69, 70, 
          71, 72, 73, 74, 75, 76, 77, 78, 79, 80, 81, 82, 83, 84, 85, 86, 87, 
          88, 89, 90, 91, 92, 93, 94, 95, 96, 97, 98, 99, 100, 101, 102, 103, 
          104, 105, 106, 107, 108, 109, 110, 111, 112, 113, 114, 115, 116, 
          117, 118, 119, 120, 121, 122, 123, 124, 125, 126, 127, 128, 129, 
          130, 131, 132, 133, 134, 135, 136, 137, 138, 139, 140, 141, 142, 
          143, 144, 145, 146, 147, 148, 149, 150, 151, 152, 153, 154, 155, 
          156, 157, 158, 159, 160, 161, 162, 163, 164, 165, 166, 167, 168 ], 
      [ 169, 170, 171, 172, 173, 174, 175, 176, 177, 178, 179, 180, 181, 182, 
          183, 184, 185, 186, 187, 188, 189, 190, 191, 192, 193, 194, 195, 
          196 ], 
      [ 197, 198, 199, 200, 201, 202, 203, 204, 205, 206, 207, 208, 209, 210, 
          211, 212, 213, 214, 215, 216, 217, 218, 219, 220, 221, 222, 223, 
          224, 225, 226, 227, 228, 229, 230, 231, 232, 233, 234, 235, 236, 
          237, 238, 239, 240, 241, 242, 243, 244, 245, 246, 247, 248, 249, 
          250, 251, 252, 253, 254, 255, 256, 257, 258, 259, 260 ] ] )
gap> r1:=SearchPFilterRowBlocks(bpQ,[1..92],[1..168],3);;
gap> Length(r1);
4374
gap> r2:=SearchPFilterRowBlocks(bpQ,[93..110],[169..196],3,[-1..1]);;
gap> Length(r2);
3936
gap> r3:=SearchPFilterRowBlocks(bpQ,[111..145],[197..260],3);;
gap> Length(r3);
1224
gap> r13:=SearchPMergeRowBlock(r1{[1..300]},r3{[1..300]});;
gap> Length(r13);
24
gap> r132:=SearchPMergeRowBlock(r13,r2);;
gap> P1:=r132{[1..168]};;
gap> P2:=r132{[233..260]};;
gap> P3:=r132{[169..232]};;
gap> P:=Concatenation([P1,P2,P3]);;
gap> Size(P); # 260x260 matrix
260
gap> Determinant(P);
1
gap> StablyPermutationFCheckMatFromBase(PSL28J[11],mi,Nlist(l)+I13[3],Plist(l)+I13[3],P);
# T is stably k-rational
true

gap> List(PSL28H,StructureDescription); # for checking core-free H
[ "1", "C2", "C3", "C2 x C2", "S3", "C7", "C2 x C2 x C2", "C9", "D14",
  "D18", "(C2 x C2 x C2) : C7" ]
gap> List([7,11],x->StructureDescription(PSL28H[x])); 
[ "C2 x C2 x C2", "(C2 x C2 x C2) : C7" ]
\end{verbatim}

\section{GAP algorithms}\label{GAPalg}

The following GAP \cite{GAP} algorithms and related ones are available as {\tt FlabbyResolutionFromBase.gap} 
in \cite{RatProbNorm1Tori}. 

\begin{verbatim}
ConjugacyClassesSubgroups2:= function(g)
    Reset(GlobalMersenneTwister);
    Reset(GlobalRandomSource);
    return ConjugacyClassesSubgroups(g);
end;

Hcandidates:= function(G)
    local Gcs,GN;
    Gcs:=ConjugacyClassesSubgroups2(G);
    GN:=NormalSubgroups(G);
    return Filtered(List(Gcs,Representative),x->Length(Filtered(GN,y->IsSubgroup(x,y)))=1);;
end;

Norm1TorusJTransitiveGroup:= function(d,n)
    local l,T,M;
    l:=Concatenation(IdentityMat(d-1),[-List([1..d-1],One)]);
    T:=TransitiveGroup(d,n);
    M:=List(GeneratorsOfGroup(T),x->List([1..d-1],y->l[y^x]));
    return Group(M);
end;

Norm1TorusJPermutationGroup:= function(G)
    local d,l,T,M;
    d:=NrMovedPoints(G);
    l:=Concatenation(IdentityMat(d-1),[-List([1..d-1],One)]);
    M:=List(GeneratorsOfGroup(G),x->List([1..d-1],y->l[y^x]));
    return Group(M);
end;

Norm1TorusJCoset:= function(g,h)
    local gg,hg,phg,d,l,M;
    gg:=GeneratorsOfGroup(g);
    hg:=SortedList(RightCosets(g,h));
    phg:=List(gg,x->Permutation(x,hg,OnRight));
    d:=Length(hg);
    l:=Concatenation(IdentityMat(d-1),[-List([1..d-1],One)]);
    M:=List(phg,x->List([1..d-1],y->l[y^x]));
    return Group(M);
end;
\end{verbatim}

%

\end{document}